# NONPARAMETRIC ESTIMATION OF CORRELATION FUNCTIONS IN LONGITUDINAL AND SPATIAL DATA, WITH APPLICATION TO COLON CARCINOGENESIS EXPERIMENTS[1]

BY YEHUA LI, NAISYIN WANG, MEEYOUNG HONG, NANCY D. TURNER,
JOANNE R. LUPTON AND RAYMOND J. CARROLL

*University of Georgia and Texas A&M University*

In longitudinal and spatial studies, observations often demonstrate strong correlations that are stationary in time or distance lags, and the times or locations of these data being sampled may not be homogeneous. We propose a nonparametric estimator of the correlation function in such data, using kernel methods. We develop a pointwise asymptotic normal distribution for the proposed estimator, when the number of subjects is fixed and the number of vectors or functions within each subject goes to infinity. Based on the asymptotic theory, we propose a weighted block bootstrapping method for making inferences about the correlation function, where the weights account for the inhomogeneity of the distribution of the times or locations. The method is applied to a data set from a colon carcinogenesis study, in which colonic crypts were sampled from a piece of colon segment from each of the 12 rats in the experiment and the expression level of p27, an important cell cycle protein, was then measured for each cell within the sampled crypts. A simulation study is also provided to illustrate the numerical performance of the proposed method.

**1. Introduction.** This paper concerns kernel-based nonparametric estimation of covariance and correlation functions. Our methods and theory are applicable to longitudinal and spatial data as well as time series data, where observations within the same subject at different time points or locations

Received October 2005; revised December 2006.
[1]Supported by the National Cancer Institute (grants CA10462, CA74552, CA57030, CA61750 and CA82907), by the National Space Biomedical Research Institute (grant NASA NCC9-58) and by the Texas A&M Center for Environmental and Rural Health via a grant from the National Institute of Environmental Health Sciences (P30-ES09106).

have strong correlations, which are stationary in time or distance lags. The structure for the observation at a particular time or location within one subject can be very general, for example, a vector or even a function.

Our study arises from a colon carcinogenesis experiment. The biomarker that we are interested in is p27, which is a life cycle protein that affects cell apoptosis, proliferation and differentiation. An important goal of the study is to understand the function of p27 in the early stage of the cancer development process. In the experiment, 12 rats were administered azoxymethane (AOM), which is a colon specific carcinogen. After 24 hours, the rats were terminated and a segment of colon tissue was excised from each rat. About 20 colonic crypts were randomly picked along a linear slice on the colon segment. The physical distances between the crypts were measured. Then, within each crypt, we measured cells at different depths within the crypts, and then the expression level of p27 was measured for each cell within the chosen crypts. In this data set, crypts are naturally functional data (Ramsay and Silverman [13]), in that the responses within a crypt are coordinated by cell depths. There is a literature about similar data, for example, Morris et al. [11].

However, in this paper we will be focused on a very different perspective. In this application the spatial correlation between crypts is of biological interest, because it helps answer the question: if we observe a crypt with high p27 expression, how likely are the neighboring crypts to have high p27 expression? We will phrase much of our discussion in terms of this example, but as seen in later sections, we have a quite general structure that includes time series as a special case. In that context, the asymptotic theory is as the number of "time series locations," that is, crypts, increases to infinity.

Although motivated by a very specific problem, nonparametric covariance/correlation estimators are worth being investigated in their own right. They can be used in a statistical analysis as: (a) an exploratory device to help formulate a parametric model, (b) an intermediate tool to do spatial prediction (kriging), (c) a diagnostic for parametric models and (d) a robust tool to test correlation. Understanding the theoretical properties of the nonparametric estimator is important under any of these situations. A limiting distribution theory would be especially valuable for purpose (d).

There is previous work on the subject of nonparametric covariance estimation. Hall, Fisher and Hoffmann [7] developed an asymptotic convergence rate of a kernel covariance estimator in a time series setting. They required not only an increasing time domain, but increasingly denser observations. Diggle and Verbyla [5] suggested a kernel-weighted local linear regression estimator for estimating the nonstationary variogram in longitudinal data, without developing asymptotic theory. Guan, Sherman and Calvin [6] used a kernel variogram estimator when assessing isotropy in geostatistics data. They proved asymptotic normality for their kernel variogram estimator in a geostatistics setting, where they required the spatial locations to be sampled from the field according to a two-dimensional homogeneous Poisson process.

As we will show below and as implied by the result from Guan, Sherman and Calvin [6], if the observations locations (or times) in the design are random, Hall's assumption, namely that the number of observations on a unit domain goes to infinity, is too restrictive and not necessary. However, in the setting of Guan, Sherman and Calvin [6], given the sample size, spatial locations are uniformly distributed within the field, which does not fit our problem, where crypt locations within a rat are, in fact, not even close to uniformly distributed.

Our paper differs from the previous work on the kernel covariance estimators in the following ways. First, our approach accommodates more complex data structure at each location or time. Second, we allow the spatial locations to be sampled in an inhomogeneous way, and as we will show below this inhomogeneity will affect the asymptotic results and inference procedures. In doing so, we generalize the setting of Guan, Sherman and Calvin [6], and link it to the setting of Hall, Fisher and Hoffmann [7]. Also, Guan, Sherman and Calvin [6] is mainly concerned with comparing variograms on a few preselected distance lags; we, on the other hand, are more interested in the correlation as a function. Third, we propose an inference procedure based upon our theory, thus filling a gap in the previous literature.

The paper is organized as follows. Section 2 introduces our model assumptions and estimators, while asymptotic results are given in Section 3. An analysis of the motivating data is given in Section 4, where we also discuss bandwidth selection and standard error estimation. Section 5 describes simulation studies, and final comments are given in Section 6. All proofs are given in the Appendix.

**2. Model assumptions and estimators.** The data considered here have the following structure.

- There are $r = 1, \ldots, R$ independent subjects, which in our example are rats.
- The data for each subject have two levels. The first level has an increasing domain, as in time series or spatial statistics, and are the crypts in our example. We label this first level as a "unit," and it is these units that have time series or spatial structure in their locations. Within each subject, there are $i = 1, \ldots, N_r$ such units.
- The second level of the data consists of observations within each of the primary units. In our case, these are the cells within the primary units, the colonic crypts. We will label this secondary level as the "subunits," which are labeled with locations. The locations with the subunits are on the interval $[0, 1]$. For simplicity, we will assume there are exactly $m$ subunits (cells) within each unit (crypt), with the $j$th subunit having location (relative cell depth) $x = (j-1)/(m-1)$. However, all theories and methods in our paper will go through if the subunits take the form of an arbitrary finite set.

- For $m = 1$, define $x$ to be fixed at 0. It is analogous to the time series setting of Hall, Fisher and Hoffman [7] or the spatial setting of Guan, Sherman and Calvin [6].

Let $\Theta(s, x)$ be a random field on $\mathcal{T} \times \mathcal{X}$, where $s$ is the unit (crypt) location and $x$ is the subunit (cell) location, so that $\mathcal{T} = [0, \infty)$, $\mathcal{X} = \{(j-1)/(m-1), j = 1, \ldots, m\}$. Assume that $\Theta_r(\cdot, \cdot)$, $r = 1, \ldots, R$, are independent realizations of $\Theta(\cdot, \cdot)$. We use the short-hand notation $\Theta_{ri}(x) = \Theta_r(S_{ri}, x)$, where $S_{ri}$ is the location of the $i$th unit (crypt) within the $r$th subject (rat). Our model for the observed data is that

$$Y_{rij} = \Theta_{ri}(x_j) + \varepsilon_{rij}, \tag{1}$$

where $Y$ is the response (logarithm of p27 level), $\varepsilon_{rij}$ are zero-mean uncorrelated measurement errors with variance $\sigma_\varepsilon^2$, $r = 1, \ldots, R$, $i = 1, \ldots, N_r$ and $j = 1, \ldots, m$ are the indices for subjects (rats), units (crypts) and subunits (cells). Define $\Psi_r(\cdot) = E_r\{\Theta_{ri}(\cdot)\}$ to be the subject-level mean, and the notation "$E_r$" refers to expectation conditional on the subject. Another way to understand $\Psi_r(\cdot)$ is to decompose the random field $\Theta_r(\cdot, \cdot)$ into the random effect model $\Theta_{ri}(x) = \Psi_r(x) + \Lambda_{ri}(x)$, where $\Psi_r(\cdot)$ is the fixed subject effect and $\Lambda_{ri}$ is the zero-mean, spatially correlated unit effect.

Within each subject, we assume that the correlation of the mean unit (crypt) level functions is stationary over the distances between the units. In addition, the covariance between unit locations $(s_1, s_2)$ at subunit (cell) locations $(x_1, x_2)$ is assumed to have the form

$$\mathcal{V}\{x_1, x_2, \Delta\} = E[\{\Theta_r(s_1, x_1) - \Psi_r(x_1)\}\{\Theta_r(s_2, x_2) - \Psi_r(x_2)\}], \tag{2}$$

where $\Delta = s_1 - s_2$. While we develop general results for model (2), in many cases it is reasonable to assume that the covariance function is separable, that is,

$$\mathcal{V}(x_1, x_2, \Delta) = G(x_1, x_2)\rho(\Delta). \tag{3}$$

When the covariance function is separable, the correlation function at the unit-level, $\rho(\cdot)$, is of interest in itself. In our application, $\rho(\cdot)$ is the correlation between crypts. We provide an estimator of $\rho(\cdot)$ as well as an asymptotic theory for that estimator.

A first estimator for the covariance function has the form

$$\widehat{\mathcal{V}}(x_j, x_l, \Delta) = \left[\sum_r \sum_i \sum_{k \neq i} K_h\{\Delta_r(i, k) - \Delta\}(Y_{rij} - \overline{Y}_{r \cdot j})(Y_{rkl} - \overline{Y}_{r \cdot l})\right]$$

(4)
$$\times \left[\sum_r \sum_i \sum_{k\neq i} K_h\{\Delta_r(i,k) - \Delta\}\right]^{-1},$$

where $\overline{Y}_{r\cdot j} = N_r^{-1}\sum_{i=1}^{N_r} Y_{rij}$, $\Delta_r(i,k) = S_{ri} - S_{rk}$ and $K_h(\cdot) = h^{-1}K(\cdot/h)$ with $K$ being a kernel function satisfying the conditions in Section 3.

It is usually reasonable to assume that $\mathcal{V}(x_1, x_2, \Delta)$ has some symmetry property, that it is an even function in $\Delta$ and $\mathcal{V}(x_1, x_2, \Delta) = \mathcal{V}(x_2, x_1, \Delta)$. However, the estimator defined in (4) does not enjoy this property. To see this, we observe that, for $x_j \neq x_l$, although $(Y_{rij} - \overline{Y}_{r\cdot j})(Y_{rkl} - \overline{Y}_{r\cdot l})$ and $(Y_{ril} - \overline{Y}_{r\cdot l})(Y_{rkj} - \overline{Y}_{r\cdot j})$ estimate the same thing, they only contribute to $\widehat{\mathcal{V}}(x_j, x_l, \Delta)$ and $\widehat{\mathcal{V}}(x_j, x_l, -\Delta)$, respectively. We also observe that $\widehat{\mathcal{V}}(x_1, x_2, \Delta) = \widehat{\mathcal{V}}(x_2, x_1, -\Delta)$.

To correct the asymmetry of the covariance estimator, for $\Delta \geq 0$, define

$$\widetilde{\mathcal{V}}(x_j, x_l, \Delta) = \left[\sum_r \sum_i \sum_{k\neq i} K_h\{|\Delta_r(i,k)| - \Delta\}(Y_{rij} - \overline{Y}_{r\cdot j})(Y_{rkl} - \overline{Y}_{r\cdot l})\right]$$
(5)
$$\times \left[\sum_r \sum_i \sum_{k\neq i} K_h\{|\Delta_r(i,k)| - \Delta\}\right]^{-1},$$

and let $\widetilde{\mathcal{V}}(x_j, x_l, \Delta) = \widetilde{\mathcal{V}}(x_j, x_l, -\Delta)$ for $\Delta < 0$. As shown in the proof of Theorem 2, for a fixed $\Delta \neq 0$, $\widetilde{\mathcal{V}}(x_1, x_2, \Delta)$ is asymptotically equivalent to $\{\widehat{\mathcal{V}}(x_1, x_2, \Delta) + \widehat{\mathcal{V}}(x_1, x_2, -\Delta)\}/2$.

In addition, when the separable structure (3) is assumed, define the estimator for the within-unit covariance as

(6) $$\widehat{G}(x_1, x_2) = \widetilde{\mathcal{V}}(x_1, x_2, 0),$$

and the estimator for the correlation function as

(7) $$\widehat{\rho}(\Delta) = \left\{\sum_{x_1 \in \mathcal{X}} \sum_{x_2 \leq x_1} \widetilde{\mathcal{V}}(x_1, x_2, \Delta)\right\} \Big/ \left\{\sum_{x_1 \in \mathcal{X}} \sum_{x_2 \leq x_1} \widehat{G}(x_1, x_2)\right\}.$$

**3. Asymptotic results.** The following are our model assumptions. Each subject (rat) is of length $L$, where in our example $L$ is the length of the segment of tissue from each rat. The units (crypts) are located on the interval $[0, L]$, and in our asymptotics we let $L \to \infty$, so that we have an increasing domain. Suppose that the positions of the units (crypts) within the $r$th subject (rat) are $S_{r1}, \ldots, S_{rN_r}$, where the $S_{ri}$'s are points from an inhomogeneous Poisson process on $[0, L]$. Then $\Delta_{r,ik} = S_{ri} - S_{rk}$. The definition of an inhomogeneous Poisson process is adopted from Cressie [3]. We assume the inhomogeneous Poisson process has a local intensity $\nu g^*(s)$, where $\nu$ is a positive constant and $g^*(s) = g(s/L)$ for a continuous density function $g(\cdot)$ on $[0, 1]$.

A special case of our setting is that $g(\cdot)$ is a uniform density function and the units (crypts) are sampled according to a homogeneous Poisson process. This is the setting investigated in Guan, Sherman and Calvin [6]. Our setting resembles that of Hall, Fisher and Hoffmann [7] in the sense that we also model the unit locations as random variables with the same distribution: in our setting, the number of units within a subject (rat) is $N_r \sim \text{Poisson}(\nu L)$, and given $N_r$, $S_{r1}/L, \ldots, S_{r,N_r}/L$ are independent and identically distributed with density $g(\cdot)$. By properties of Poisson processes, $N_r/L = O(\nu)$ almost surely, as $L \to \infty$, that is, the number of units (crypts) on a unit length tends to a constant. It is worth noting that Hall, Fisher and Hoffmann [7] required this ratio to go to infinity. We require less samples on the domain than do Hall, Fisher and Hoffmann [7].

In what follows, we provide a list of definitions and conditions.

1. We assume that $g(\cdot)$ is continuous and $c_1 \geq g(t) \geq c_2 > 0$ for all $t \in [0, 1]$. Suppose $t_i$, $i = 1, 2, 3, 4$, are independent random variables with density $g(\cdot)$, define $f_1$, $f_2$ and $f_3$ to be the densities for $t_1 - t_2$, $(t_1 - t_2, t_3 - t_2)$ and $(t_1 - t_2, t_3 - t_4, t_2 - t_4)$, respectively. Since $g(\cdot)$ is bounded, one can easily derive that $f_1(0)$, $f_2(0, 0)$ and $f_3(0, 0, 0)$ are positive. We also assume that $f_1$ and $f_2$ are Lipschitz continuous in the neighborhood of $\mathbf{0}$, that is, $|f_i(\mathbf{u}) - f_i(\mathbf{0})| \leq \lambda_i \|\mathbf{u}\|$, for $\forall \mathbf{u}$ and some fixed constants $\lambda_i > 0$, $i = 1, 2$.
2. Assume $\mathcal{V}(x_1, x_2, \Delta)$ has two bounded continuous partial derivatives in $\Delta$, and that $\sup_{x_1, x_2} \int |\mathcal{V}(x_1, x_2, \Delta)| \, d\Delta < \infty$.
3. Let
$$\mathcal{M}(x_1, x_2, x_3, x_4, u, v, w)$$
$$= E_r[\{\Theta_{ri_1}(x_1) - \Psi_r(x_1)\}\{\Theta_{ri_2}(x_2) - \Psi_r(x_2)\}\{\Theta_{ri_3}(x_3) - \Psi_r(x_3)\}$$
$$\times \{\Theta_{ri_4}(x_4) - \Psi_r(x_4)\} | \Delta_r(i_1, i_2) = u,$$
$$\Delta_r(i_3, i_4) = v, \Delta_r(i_2, i_4) = w]$$
$$- \mathcal{V}(x_1, x_2, u)\mathcal{V}(x_3, x_4, v).$$

We assume $\mathcal{M}$ has bounded partial derivatives in $u$, $v$ and $w$, and

$$(8) \qquad \sup_{x_1, x_2, x_3, x_4, u, v} \int |\mathcal{M}(x_1, x_2, x_3, x_4, u, v, w)| \, dw < \infty.$$

4. Denote $b_r(x_1, x_2, \Delta) = L^{-1} \sum_i \sum_{k \neq i} K_h\{\Delta - \Delta_r(i, k)\}\{Y_r(S_{ri}, x_1) - \Psi_r(x_1)\}\{Y_r(S_{rk}, x_2) - \Psi_r(x_2)\}$. We assume that, for any fixed $\Delta$, for some $\eta > 0$,

$$(9) \qquad \sup_{L, x_1, x_2} E(|\operatorname{var}^{-1/2}\{b_r(x_1, x_2, \Delta)\}[b_r(x_1, x_2, \Delta) - E\{b_r(x_1, x_2, \Delta)\}]|^{2+\eta})$$
$$\leq C_\eta < \infty.$$

5. Let $\mathcal{F}(T)$ be the $\sigma$-algebra generated by $\{\Theta(s, x), s \in T, x \in \mathcal{X}\}$, for any Borel set $T \subset \mathcal{T}$. Assume that the random field satisfies the mixing condition

$$\alpha(\tau) = \sup_t[|P(A_1 \cap A_2) - P(A_1)P(A_2)| : A_1 \in \mathcal{F}\{[0, t]\},$$

(10) $$A_2 \in \mathcal{F}\{[t + \tau, \infty)\}]$$

$$= O(\tau^{-\delta}) \quad \text{for some } \delta > 0.$$

6. The kernel function $K$ is a symmetric, continuous probability density function supported on $[-1, 1]$. Define $\sigma_K^2 = \int u^2 K(u) \, du$ and $R_K = \int K^2(v) \, dv$.
7. Assume that $m$ and $R$ are fixed numbers, $L \to \infty$, $h \to 0$, $Lh \to \infty$ and $Lh^5 = O(1)$.

In assumption 1, we are imposing some regularity conditions on $g$ and $f_i$. In fact, when $g$ is differentiable $f_i$ are piecewise differentiable, but usually not differentiable at $\mathbf{0}$. However, the Lipschitz conditions on $f_1$ and $f_2$ are easily satisfied when, for example, $g$ is Lipschitz continuous on $[0, 1]$.

Since we are estimating the covariance function, which is the second moment function, we need a regularity condition on the fourth moment function as in (8). Condition (9) may seem strong at the first sight, but it is simply a condition that bounds the tail probability of our statistics. For example, if we have an assumption analogous to (8) for the eighth moment of $\Theta_r(s, x)$, we can use arguments as in Lemma A.3 to show that $E([b_r(x_1, x_2, \Delta) - E\{b_r(x_1, x_2, \Delta)\}]^4) = O(L^{-3}h^{-3})$, and therefore condition (9) is satisfied for $\eta = 2$. In general, when the distribution of $\Theta$ is neither too skewed nor has a much heavier tail than that of the Gaussian, equation (9) will be satisfied. Assumption 6 and 7 are standard in the literature of kernel estimators.

Denote $\mathcal{V}^{(0,0,2)}(x_1, x_2, \Delta) = \partial^2 \mathcal{V}(x_1, x_2, \Delta)/\partial \Delta^2$. Let $\mathcal{V}(\Delta)$, $\widehat{\mathcal{V}}(\Delta)$ and $\widetilde{\mathcal{V}}(\Delta)$ denote the vectors collecting $\mathcal{V}(x_1, x_2, \Delta)$, $\widehat{\mathcal{V}}(x_1, x_2, \Delta)$ and $\widetilde{\mathcal{V}}(x_1, x_2, \Delta)$, respectively, for all distinct pairs of $(x_1, x_2)$. The following are our main theoretical results. All proofs are provided in the Appendix. Note that Theorem 1 refers to $\widehat{\mathcal{V}}(\cdot)$ in (4), while Theorem 2 refers to $\widetilde{\mathcal{V}}(\cdot)$ in (5).

THEOREM 1. *Under assumptions 1–7, for $\Delta \neq \Delta'$, we have*

$$(RLh)^{1/2} \begin{bmatrix} \widehat{\mathcal{V}}(\Delta) - \mathcal{V}(\Delta) - \text{bias}\{\widehat{\mathcal{V}}(\Delta)\} \\ \widehat{\mathcal{V}}(\Delta') - \mathcal{V}(\Delta') - \text{bias}\{\widehat{\mathcal{V}}(\Delta')\} \end{bmatrix}$$

$$\Rightarrow \text{Normal}\left[0, \{\nu^2 f_1(0)\}^{-1} \begin{pmatrix} \Sigma(\Delta) & C(\Delta, \Delta') \\ C^T(\Delta, \Delta') & \Sigma(\Delta') \end{pmatrix}\right],$$

where the asymptotic bias $\mathrm{bias}\{\widehat{\mathcal{V}}(\Delta)\}$ is a vector having entries $\mathrm{bias}\{\widehat{\mathcal{V}}(x_1, x_2, \Delta)\} = \sigma_K^2 \mathcal{V}^{(0,0,2)}(x_1, x_2, \Delta) h^2/2$, $\Sigma(\Delta)$ is the covariance matrix with the entry corresponding to $\mathrm{cov}\{\widehat{\mathcal{V}}(x_1, x_2, \Delta), \widehat{\mathcal{V}}(x_3, x_4, \Delta)\}$ equal to $R_K\{\mathcal{M}(x_1, x_2, x_3, x_4, \Delta, \Delta, 0) + I(x_2 = x_4)\sigma_\varepsilon^2 \mathcal{V}(x_1, x_3, 0) + I(x_1 = x_3)\sigma_\varepsilon^2 \mathcal{V}(x_2, x_4, 0) + I(x_1 = x_3, x_2 = x_4)\sigma_\varepsilon^4\} + I(\Delta = 0) R_K \{\mathcal{M}(x_1, x_2, x_3, x_4, 0, 0, 0) + I(x_1 = x_4)\sigma_\varepsilon^2 \mathcal{V}(x_2, x_3, 0) + I(x_2 = x_3)\sigma_\varepsilon^2 \mathcal{V}(x_1, x_4, 0) + I(x_1 = x_4, x_2 = x_3)\sigma_\varepsilon^4\}$ and $C(\Delta, \Delta')$ is the matrix with the entry corresponding to $\mathrm{cov}\{\widehat{\mathcal{V}}(x_1, x_2, \Delta), \widehat{\mathcal{V}}(x_3, x_4, \Delta')\}$ equal to $I(\Delta' = -\Delta)\{\mathcal{M}(x_1, x_2, x_3, x_4, \Delta, -\Delta, -\Delta) + I(x_2 = x_3)\sigma_\varepsilon^2 \mathcal{V}(x_1, x_4, 0) + I(x_1 = x_4) \times \sigma_\varepsilon^2 \mathcal{V}(x_2, x_3, 0) + I(x_1 = x_4, x_2 = x_3)\sigma_\varepsilon^4\}$.

THEOREM 2. *Under assumptions 1–7, for $\Delta \neq \pm \Delta'$, we have*

$$(RLh)^{1/2} \begin{bmatrix} \widetilde{\mathcal{V}}(\Delta) - \mathcal{V}(\Delta) - \mathrm{bias}\{\widetilde{\mathcal{V}}(\Delta)\} \\ \widetilde{\mathcal{V}}(\Delta') - \mathcal{V}(\Delta') - \mathrm{bias}\{\widetilde{\mathcal{V}}(\Delta')\} \end{bmatrix}$$
$$\Rightarrow \mathrm{Normal}\left[0, \{\nu^2 f_1(0)\}^{-1} \begin{pmatrix} \Omega(\Delta) & 0 \\ 0 & \Omega(\Delta') \end{pmatrix}\right],$$

where $\mathrm{bias}\{\widetilde{\mathcal{V}}(\Delta)\}$ is a vector with entries $\mathrm{bias}\{\widetilde{\mathcal{V}}(x_1, x_2, \Delta)\} = \sigma_K^2 \mathcal{V}^{(0,0,2)}(x_1, x_2, \Delta) h^2/2$, $\Omega(\Delta)$ is the covariance matrix with the entry corresponding to $\mathrm{cov}\{\widetilde{\mathcal{V}}(x_1, x_2, \Delta), \widetilde{\mathcal{V}}(x_3, x_4, \Delta)\}$ equal to $(1/2) R_K\{\mathcal{M}(x_1, x_2, x_3, x_4, \Delta, \Delta, 0) + \mathcal{M}(x_1, x_2, x_3, x_4, \Delta, -\Delta, -\Delta) + I(x_2 = x_4)\sigma_\varepsilon^2 \mathcal{V}(x_1, x_3, 0) + I(x_1 = x_3)\sigma_\varepsilon^2 \mathcal{V}(x_2, x_4, 0) + I(x_1 = x_3, x_2 = x_4)\sigma_\varepsilon^4 + I(x_2 = x_3)\sigma_\varepsilon^2 \mathcal{V}(x_1, x_4, 0) + I(x_1 = x_4)\sigma_\varepsilon^2 \mathcal{V}(x_2, x_3, 0) + I(x_1 = x_4, x_2 = x_3)\sigma_\varepsilon^4\} + I(\Delta = 0)(1/2) R_K \{2\mathcal{M}(x_1, x_2, x_3, x_4, 0, 0, 0) + I(x_2 = x_4)\sigma_\varepsilon^2 \mathcal{V}(x_1, x_3, 0) + I(x_1 = x_3)\sigma_\varepsilon^2 \mathcal{V}(x_2, x_4, 0) + I(x_1 = x_3, x_2 = x_4)\sigma_\varepsilon^4 + I(x_2 = x_3)\sigma_\varepsilon^2 \mathcal{V}(x_1, x_4, 0) + I(x_1 = x_4)\sigma_\varepsilon^2 \mathcal{V}(x_2, x_3, 0) + I(x_1 = x_4, x_2 = x_3)\sigma_\varepsilon^4\}$.

COROLLARY 1. *Suppose the covariance function has the separable structure in (3) with $\sum_{x_1} \sum_{x_2 \leq x_1} G(x_1, x_2) \neq 0$ and $\widehat{\rho}(\Delta)$ defined in (7). Then for $\Delta \neq 0$, we have*

$$(RLh)^{1/2}[\widehat{\rho}(\Delta) - \rho(\Delta) - \mathrm{bias}\{\widehat{\rho}(\Delta)\}] \Rightarrow \mathrm{Normal}[0, \{\nu^2 f_1(0)\}^{-1} \sigma_\rho^2(\Delta)],$$

where $\mathrm{bias}\{\widehat{\rho}(\Delta)\} = \{\rho^{(2)}(\Delta) - \rho(\Delta)\rho^{(2)}(0)\}\sigma_K^2 h^2/2$ is the asymptotic bias of $\widehat{\rho}(\Delta)$ and $\sigma_\rho^2(\Delta) = \{\sum_{x_1} \sum_{x_2 \leq x_1} G(x_1, x_2)\}^{-2}\{\mathbf{1}^T \Omega(\Delta) \mathbf{1} + \rho^2(\Delta) \mathbf{1}^T \Omega(0) \mathbf{1}\}$.

REMARK. The measurement errors in (1) affect the covariance estimator mainly though the nugget effect [3]. In our covariance estimators (4) and (5), we eliminate the nugget effect by excluding the $k = i$ terms in the summation. As a result, the measurement errors do not introduce bias to our covariance estimators. However, they do affect the variation of the covariance estimators and hence the correlation estimator, as seen by the fact that $\sigma_\varepsilon^2$ is in the variance expressions for all our estimators.

**4. Data analysis.** In this section we apply our methods to study the between-crypt dependence in the carcinogenesis experiment. Recall that the main subjects are rats, the units of interest are colonic crypts and the subunits within a unit are cells at which we observe the logarithms of p27 in a cell. The subunit locations that we work with in this illustration are at $x = 0, 0.1, 0.2, \ldots, 1.0$. We discuss three key issues in our analysis, namely bandwidth selection, standard error estimation and positive semidefinite adjustment, in the following three subsections.

4.1. *Bandwidth selection.*

4.1.1. *Global bandwidth.* Diggle and Verbyla [5] suggested a cross-validation procedure to choose the bandwidth for a kernel variogram estimator. We modify their procedure into the following two types of "leave-one-subject-out" cross-validation criteria. The first is based on prediction error without assuming any specific covariance structure and is given as

$$(11) \quad CV_1(h) = \sum_r \sum_{|\Delta_r(i,k)|<\Delta_0} \sum_{j=1}^m \sum_{l=1}^m [v_{r,ik}(x_j, x_l) - \tilde{\mathcal{V}}_{(-r)}\{x_j, x_l, \Delta_r(i,k)\}]^2,$$

where $v_{r,ik}(x_j, x_l) = (Y_{rij} - \overline{Y}_{r \cdot j})(Y_{rkl} - \overline{Y}_{r \cdot l})$, $\tilde{\mathcal{V}}_{(-r)}(x_1, x_2, \Delta)$ is the kernel covariance estimator using bandwidth $h$, as defined in (5), with all information on the $r$th subject (rat) left out. Here we focus on the range $|\Delta_r(i,k)| < \Delta_0$, where $\Delta_0$ is a prechosen cut-off point. The criterion $CV_1(h)$ thus evaluates the prediction error for different $h$ within the range of $|\Delta_r(i,k)| < \Delta_0$.

The cross-validation criterion (11) assumes no specific covariance structure, while our second cross-validation criterion takes into account the separable structure in (3) and is given as

$$(12) \quad CV_2(h) = \sum_r \sum_{|\Delta_r(i,k)|<\Delta_0} \sum_{j=1}^m \sum_{l=1}^m [v_{r,ik}(x_j, x_l) - \widehat{G}_{(-r)}(x_j, x_l)\widehat{\rho}_{(-r)}\{\Delta_r(i,k)\}]^2,$$

where $\widehat{G}_{(-r)}(x_1, x_2)$ and $\widehat{\rho}_{(-r)}(\Delta)$ are the estimators of $G$ and $\rho$ defined in (6) and (7), with the $r$th subject (rat) left out.

We evaluated both criteria to estimate the bandwidth $h$. We chose $\Delta_0 = 500$ microns. The first two columns of Table 1 give the minimum points and minimum values of the two cross-validation criteria.

By observing Table 1, we find the two criteria gave almost identical minimum values. Since the cross-validation scores are estimates of the prediction errors, the two cross-validation criteria represent prediction errors with and without the separable structure (3). The phenomenon that $CV_1(\cdot)$ and $CV_2(\cdot)$ have almost the same minimum values suggests that the separability assumption (3) fits the data well.

TABLE 1
*Outcomes of two cross-validation procedures on the carcinogenesis p27 data*

|  | Optimal $h$ | Min CV score | Min score, 2 par |
|---|---|---|---|
| $CV_1$ | 124.2334 | 6.5073 | 6.4867 |
| $CV_2$ | 122.7202 | 6.4955 | 6.4788 |

The data used in the validation are those with $\Delta$ values less than $\Delta_0 = 500$ microns. The first column gives the optimal global bandwidth, the second column gives the value of the cross-validation function at the optimal global bandwidth and the third column gives the minimum value of the cross-validation functions using two different smoothing parameters.

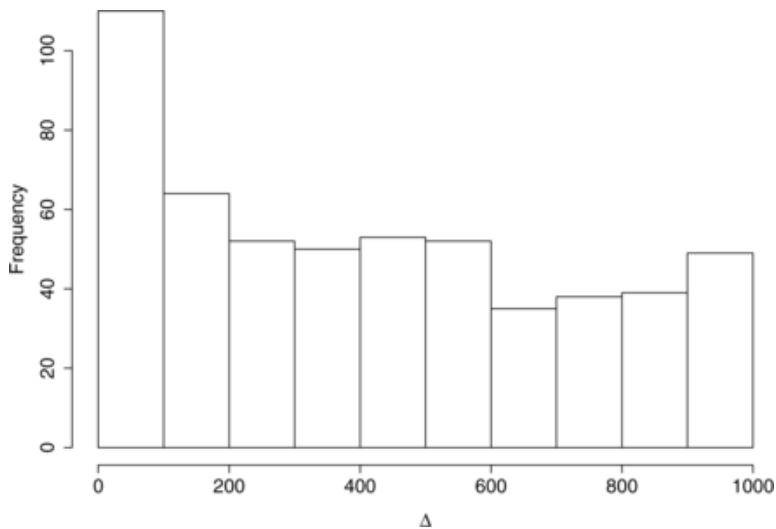

Fig. 1. *Histogram of $|\Delta_r(i,k)|$ in the carcinogenesis p27 data. $|\Delta|$ less than 1000 microns are considered.*

4.1.2. *Two bandwidths.* The independent variables in the kernel estimator are $|\Delta_r(i,k)|$ for all pairs of crypts within one subject. As shown in Figure 1, the distribution of $|\Delta_r(i,k)|$ that are less than 1000 microns, even more than the target range of interest, is locally somewhat akin to a uniform distribution.

As a robustness check on the global bandwidth, we repeated our analysis, except we used one bandwidth for $|\Delta| \leq 200$ microns, and we used a second bandwidth for $|\Delta| > 200$, and then repeated the cross-validation calculations in (11) and (12). The minimum values of the two cross-validation criteria are reported in the third column of Table 1.

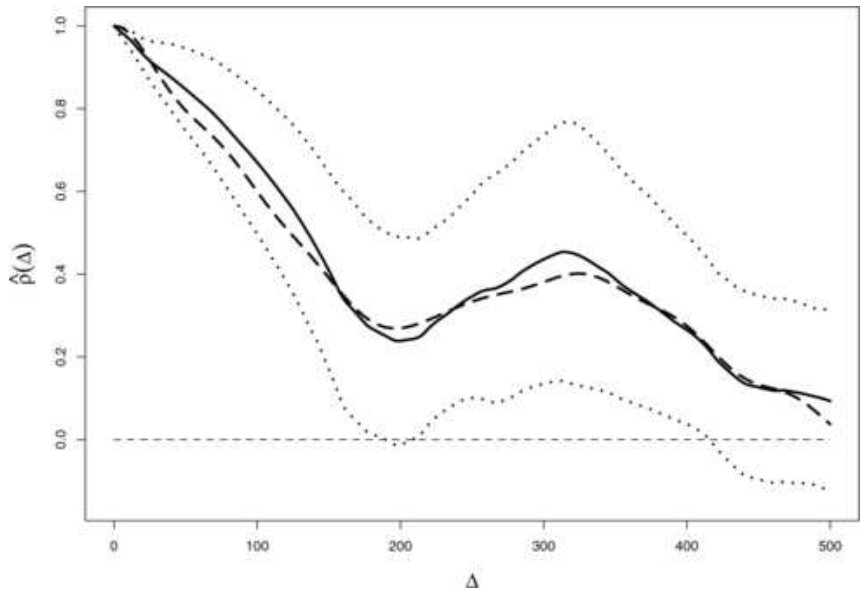

Fig. 2. *Estimated correlation function for the carcinogenesis p27 data. The solid curve is $\widehat{\rho}(\Delta)$, the dotted curves are $\widehat{\rho}(\Delta) \pm 2\widehat{SD}\{\widehat{\rho}(\Delta)\}$, and the dashed curve is the positive semidefinite adjusted estimate, $\widetilde{\rho}(\Delta)$.*

Comparing the results in columns 2 and 3 in Table 1, we find the minimum values of the cross-validation functions did not change much, that is, an extra smoothing parameter did not substantially reduce the prediction error for the domain $|\Delta| \leq 500$ microns. In other words, it appears sufficient to use a global bandwidth to estimate $\rho(\Delta)$ for $|\Delta| \leq 500$. For the following analysis, we use the bandwidth $h = 122$ microns, as suggested by $CV_2$, and the resulting estimate $\widehat{\rho}$ is shown as the solid curve in Figure 2.

4.2. *Standard error estimation.* Our primary goal in this section is to construct a standard error estimate for $\widehat{\rho}(\Delta)$.

The asymptotic variance of $\widehat{\rho}(\Delta)$ has a very complicated form, which involves the fourth moment function of the random field, $\mathcal{M}(x_1, x_2, x_3, x_4, u, v, w)$. With so many estimates of higher-order moments involved, a plug-in method, while feasible, is not desirable. We instead use a bootstrap method to estimate the variance directly.

In our model assumptions, the number of subjects (rats) $R$ is fixed, which means that bootstrapping solely on the subject level will not give a consistent estimator of the variance. Consequently, we decided to subsample within each subject. When the data are dependent, block bootstrap methods have been investigated and used, see Shao and Tu [15]). Politis and Sherman [12] also justified using a block sub-sampling method to estimate the variance of a statistic when the data are from a marked point process. Our data can be viewed as a marked inhomogeneous Poisson process. However, because of the inhomogeneity, we need to modify their procedure: when we subsample a block from each subject and compute the statistic $\widehat{\rho}(\Delta)$ by combining these blocks, the variance of the statistic depends on the corresponding local intensity at the location where each block is sampled.

By letting $R = 1$ in Corollary 1, our theory implies that if the number of units goes to infinity, each subject will provide a consistent estimator of $\rho(\Delta)$. Now, suppose the Poisson process for each subject has a different local intensity, $\nu_r g_r^*(s)$, $r = 1, \ldots, R$. With a slight modification of our theoretical derivations, one can show that

$$\left\{\sum_{r=1}^R \nu_r^2 f_{r,1}(0) L h\right\}^{1/2} [\widehat{\rho}(\Delta) - \rho(\Delta) - \text{bias}\{\widehat{\rho}(\Delta)\}] \Rightarrow \text{Normal}\{0, \sigma_\rho^2(\Delta)\},$$

where $f_{r,1}(t) = \int g_r(t+u) g_r(u) \, du$, $r = 1, \ldots, R$, are the counterparts of $f_1(t)$ used in Theorem 1, 2 and Corollary 1.

Define $A(\Delta) = \sum_r \sum_i \sum_{k \neq i} K_h\{\Delta_r(i,k) - \Delta\}$. Then by Lemma A.2,

$$A(\Delta) \Big/ \left\{\sum_{r=1}^R \nu_r^2 f_{r,1}(0) L\right\} \to 1, \qquad \text{in } L^2.$$

Note that $A(\Delta)$ here is defined slightly different from $a(\Delta)$ in Lemma A.2.

This equation suggests a natural way to construct correction weights in a weighted block bootstrap procedure that accommodates the inhomogeneous local intensity. We now propose our weighted bootstrap procedure:

1. Resample $R$ subjects (rats) with replacement.
2. Within each resampled subject, randomly subsample a block with length $L^*$.
3. Combine the $R$ blocks as our resampled data, compute $\widehat{\rho}(\Delta)$ and $A(\Delta)$ using the resampled data, with the same bandwidth $h$ as for the kernel estimator (7).
4. Repeat steps 1–3 $B$ times, denoting the results from the $b$th iteration as $\widehat{\rho}_b^*(\Delta)$ and $A_b^*(\Delta)$.
5. Obtain the estimator of the standard deviation as

$$\widehat{\text{sd}}\{\widehat{\rho}(\Delta)\} = \left[A^{-1}(\Delta) B^{-1} \sum_{b=1}^B A_b^*(\Delta) \{\widehat{\rho}_b^*(\Delta) - \widehat{\rho}_\cdot^*(\Delta)\}^2\right]^{1/2},$$

where $\widehat{\rho}_\cdot^*(\Delta) = B^{-1} \sum_{b=1}^B \widehat{\rho}_b^*(\Delta)$.

The block length $L^*$ should increase slowly with $L$. Politis and Sherman [12] proposed a block size selection procedure for dependent data on irregularly-spaced observation points which are from a homogeneous point process.

Their procedure is built on asymptotic theory and needs a good pilot block size. The good performance of the procedure often requires a fairly large sample size. The implementation could be computationally intense.

One operational idea for a moderate sample size in our context is to choose $L^*$ such that the correlation dies out outside the block but there are still a relatively large number of blocks. For this data set, we adopted a block size such that there are at least a couple of nonoverlapping blocks within each subject, and there are totally 24 nonoverlapping blocks by pooling all subjects together. In our analysis we took $L^* = 1$ cm ($=10{,}000$ microns). We also tried $L^* = 8$ mm and $L^* = 1.1$ cm; the results are very similar. We investigate the numerical performance of this simple procedure in both of the two simulation studies in Section 5, and we find it works pretty well.

The two dotted curves in Figure 2 show $\widehat{\rho} \pm 2$ standard deviation. The plot implies that the correlation is practically zero when the crypt distance exceeds 500 microns.

4.3. *Positive semidefinite adjustment.* By definition, $\rho(\Delta)$ is a stationary correlation function and therefore is positive semidefinite, that is, $\iint \rho(\Delta_1 - \Delta_2)\omega(\Delta_1) \times \omega(\Delta_2)\, d\Delta_1\, d\Delta_2 \geq 0$ for all integrable functions $\omega(\cdot)$. By Bochner's theorem, the positive semidefiniteness is equivalent to nonnegativity of the Fourier transformation of $\rho$, that is, $\rho^+(\theta) \geq 0$ for all $\theta$, where $\rho^+(\theta) = \int_{-\infty}^{\infty} \rho(\Delta)\exp(i\theta\Delta)\, d\Delta = 2\int_0^{\infty} \rho(\Delta)\cos(\theta\Delta)\, d\Delta$.

To make $\widehat{\rho}$ a valid correlation function, we apply an adjustment procedure suggested by Hall and Patil [8]. First we compute the Fourier transformation of $\widehat{\rho}(\cdot)$,

$$\widehat{\rho}^+(\theta) = 2\int_0^{\infty} \widehat{\rho}(\Delta)\cos(\theta\Delta)\, d\Delta.$$

In practice, we cannot accurately estimate $\rho(\Delta)$ for a large $\Delta$ because of data constraints. So, what we should do is multiply $\widehat{\rho}$ by a weight function $w(\Delta) \leq 1$ and let

$$\widehat{\rho}^+(\theta) = 2\int_0^{\infty} \widehat{\rho}(\Delta)w(\Delta)\cos(\theta\Delta)\, d\Delta.$$

Possible choices of $w(\cdot)$ suggested by Hall, Fisher and Hoffmann [7] are $w_1(\Delta) = I(|\Delta| \leq D)$ for some threshold value $D > 0$; and $w_2(\Delta) = 1$ if $|\Delta| < D_1$, $(D_2 - |\Delta|)/(D_2 - D_1)$ if $D_1 \leq |\Delta| \leq D_2$ and 0 if $|\Delta| > D_2$.

The next step is to make $\widehat{\rho}^+$ nonnegative and then take an inverse Fourier transformation. So, the adjusted estimator is defined as

$$\widetilde{\rho}(\Delta) = (2\pi)^{-1}\int \max\{\widehat{\rho}^+(\theta), 0\}\cos(\theta\Delta)\, d\theta.$$

The adjusted estimate of the correlation function for the colon carcinogenesis p27 data is given as the dashed curve in Figure 2.

**5. Simulation studies.** We present three simulation studies to illustrate the numerical performance of the kernel correlation estimator under different settings.

5.1. *Simulation* 1. Our first simulation study is to mimic the colon carcinogenesis data, so that the result can be inferred to evaluate the performance of our estimators in the data analysis and to justify our choice of tuning parameters.

The simulated data arise from the model

$$Y_r^*(s_{ri}, x_j) = \Theta_r^*(s_{ri}, x_j) + \varepsilon_{rij}^*,$$

where $\Theta_r^*(s, x)$ is the $r$th replicate of a zero-mean Gaussian random field $\Theta^*(s, x)$, $r = 1, \ldots, 12$. As in our data analysis, $x$ takes values in $\{0.0, 0.1, \ldots, 0.9, 1.0\}$. We used the actual unit (crypt) locations from the data as the sample locations $s_{ri}$ in the simulated data. In addition, $\Theta^*(s, x)$ has covariance structure (2) and (3), with

$$G^*(x_1, x_2) = \left(\sum_{r=1}^{12} N_r\right)^{-1} \sum_{r=1}^{12} \sum_{i=1}^{N_r} \{Y_{ri}(x_1) - \overline{Y}_{r\cdot}(x_1)\}\{Y_{ri}(x_2) - \overline{Y}_{r\cdot}(x_2)\},$$

(13)

which is computed from the data, and $\rho^*(\Delta)$ chosen from the Matérn correlation family $\rho^*(\Delta; \phi, \kappa) = \{2^{\kappa-1}\Gamma(\kappa)\}^{-1}(\Delta/\phi)^\kappa K_\kappa(\Delta/\phi)$, where $K_\kappa(\cdot)$ is the modified Bessel function; see Stein [16]. In our simulation we chose $\kappa = 1.5$ and $\phi = 120$ microns. In addition, the $\varepsilon_{rij}^*$ are independent identically distributed as Normal$(0, \sigma_{\varepsilon^*}^2)$. For $\sigma_{\varepsilon^*}^2$ we use an estimate of $\sigma_\varepsilon^2$ from the data, $\sigma_{\varepsilon^*}^2 = \frac{1}{11}\sum_{j=1}^{11}\{G^*(x_j, x_j) - \widehat{G}(x_j, x_j)\}$, where $x_j = (j-1)/10$, $j = 1, \ldots, 11$, and $G^*$ and $\widehat{G}$ are defined in (13) and (6), respectively.

For each simulated data set, we computed $\widehat{\rho}(\Delta)$ and the standard deviation estimator $\widehat{SD}\{\widehat{\rho}(\Delta)\}$ that we proposed in Section 4.2, for bandwidths $h = 120$ and 200 microns. When doing the bootstrap, we used block size $L^* = 1$ cm, as we did in the p27 data analysis. We repeated the simulation 200 times.

Figure 3 shows the means and 5% and 95% pointwise percentiles of $\widehat{\rho}$ for the two bandwidths, and compares them to the truth $\rho^*$. Obviously, as expected from the theory, the larger bandwidth incurs the bigger bias. By the plots, it seems that when $h = 120$ the kernel estimator $\widehat{\rho}$ behaves quite well. We compare the true bias from the simulation study to the asymptotic bias computed with the true correlation function $\rho^*$, under bandwidth $h = 120$. We find the differences between the two are less than 0.04. This means the bias shown in Figure 3 is explainable by our asymptotic theory.

In Figure 4, we show the pointwise standard deviation of $\widehat{\rho}$ from the simulation and the mean of the bootstrap standard deviation estimates.

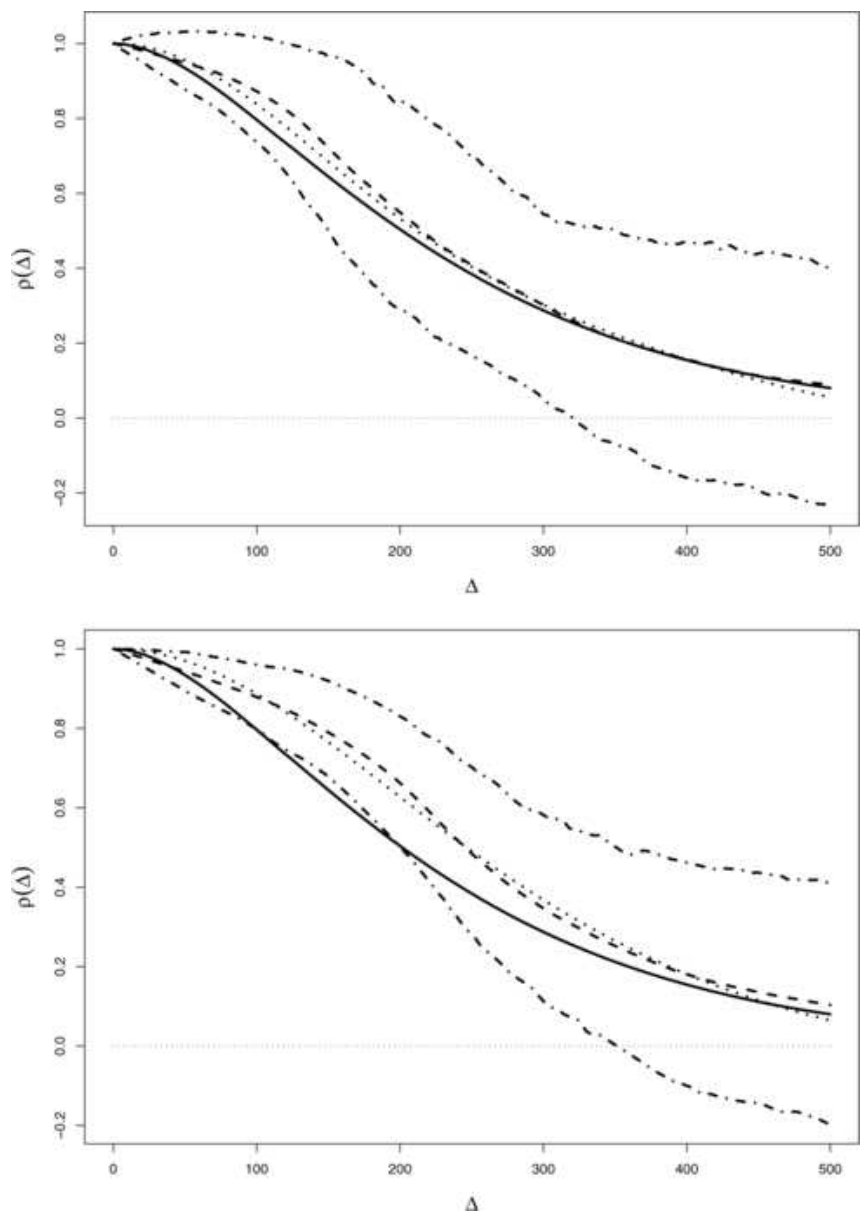

FIG. 3. *Plots of the correlation estimators in Simulation 1. Upper panel: $h = 120$; lower panel: $h = 200$. In each plot, the solid curve is the true correlation function $\rho(\cdot)$, the dashed curve is the mean of $\widehat{\rho}(\cdot)$, the dotted curve is the mean of $\widetilde{\rho}(\cdot)$, and the dot-dash curves are the 5% and 95% pointwise percentiles of $\widehat{\rho}$. $h = 200$ oversmooths the curve, hence incurs larger bias.*

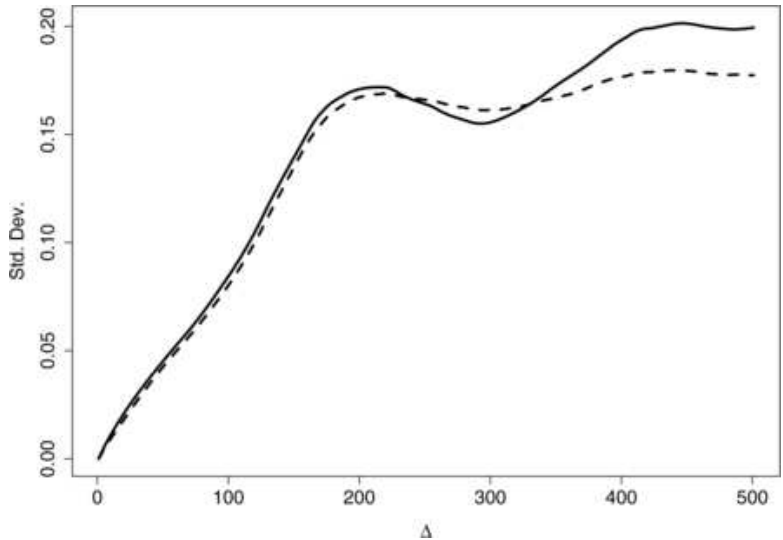

FIG. 4. *Standard deviation of $\widehat{\rho}$ in Simulation 1. The solid curve is the pointwise standard deviation of $\widehat{\rho}$ from the simulation in Section 5.1, and the dashed curve is the mean of the 200 bootstrap standard deviation estimates. The bandwidth $h = 120$ was used.*

The closeness of the two curves implies that our bootstrap procedure in Section 4.2 gives an approximately unbiased estimator of the true standard deviation, which also implies that our choice of block length, $L^* = 1$ cm, is reasonable. In our simulation, we also tried $L^* = 8$ mm and $L^* = 1.1$ cm; the results are very similar.

We applied the positive semidefinite adjustment procedure in Section 4.3 to the simulation, and the pointwise mean of $\widetilde{\rho}(\cdot)$ is also shown in Figure 3. We computed the integrated mean squared errors (IMSE) of $\widehat{\rho}$ and $\widetilde{\rho}$ up to $\Delta = 500$, and found that $\text{IMSE}(\widehat{\rho}) = 12.56$ whereas $\text{IMSE}(\widetilde{\rho}) = 8.50$. This result agrees with the theory in Hall and Patil [8] that positive semidefinite adjustment can actually improve the integrated mean squared error of the raw kernel estimator. We found that most of the improvements come from the regions where $\rho(\cdot)$ is close to 0 or 1, the areas where the procedure corrects the shape of $\widehat{\rho}$ the most due to the enforcement of positive semidefiniteness.

5.2. *Simulation 2.* As suggested by the referees, we provide a second simulation study to evaluate the finite sample numerical performance of our correlation estimator when the locations or times are from an inhomogeneous Poisson process as assumed in Section 3. Also, we choose a correlation function which is similar in shape to that obtained from the p27 data example, but is even less monotone; this clearly illustrates a situation that an "off-the-shelf" parametric model fails to fit the data. The true correlation function is given by the solid curve in the middle panel in Figure 5, while the corresponding spectral density is given in the upper panel of the same figure.

We kept the same simulation setup as those in Simulation 1, except that the spatial correlation $\rho(\Delta)$ was set to be the one given in Figure 5, the locations were sampled from an inhomogeneous Poisson process as given in Section 3 with $g$ a truncated normal density function on $[0,1]$, and we simulated only one subject on a prolonged domain $[0, L]$, with $L = 50{,}000$ and the expected number of units equal to 500. We let the bandwidth $h = 35$ and block size $L^* = 6000$ for the weighted block bootstrap procedure.

We repeated the process described above 200 times, and computed the proposed correlation estimator for each simulated data set. In the middle panel of Figure 5, the mean of our kernel correlation estimator is given by the dashed curve, while the dotted curve is the best approximation to the true correlation function from the Matérn family. As one can see, our nonparametric method can consistently estimate a nonmonotone correlation function.

In Figure 5 we also compare the mean of our bootstrap standard deviation estimator with the true pointwise standard deviation curve. We found that the proposed standard deviation estimator also works quite well given the finite sample size.

5.3. *Simulation* 3. Our third simulation study has the same setting as Simulation 2, except that the correlation function is replaced by

$$\rho(\Delta) = \frac{1}{2} \frac{\cos(\Delta/60)}{1 + |\Delta|/100} + \frac{1}{2} \exp\left(-\frac{|\Delta|}{800}\right),$$

as suggested by one referee. This correlation is given by the solid curve in Figure 6. We use this simulation to illustrate the performance of the kernel correlation estimator and its adjusted version in the case that the true correlation function is not smooth at 0.

Because this function decays to 0 with a slower rate, we present the estimate up to $\Delta = 1000$. The dashed curve in Figure 6 gives the mean of $\widehat{\rho}$ over 200 simulations, the two dotted curves give the pointwise 5% and 95% percentiles of $\widehat{\rho}$ and the dot–dash curve gives the mean of $\widetilde{\rho}$. One can see that $\widehat{\rho}$ still behaves well even though the true function is not differentiable at 0. The adjustment procedure introduced some bias, but reduced the variation. We compared the IMSE of the two estimators over $[0, 50]$ and $[0, 500]$: the IMSE values for $\widehat{\rho}$ over the two ranges are 0.40 and 6.59, while the corresponding IMSE values for $\widetilde{\rho}$ are 0.19 and 4.53. The adjustment procedure improved the IMSE on both ranges, even in an area close to the origin.

This simulation study shows that our estimators work even when the differentiability assumption in Section 3 is mildly violated.

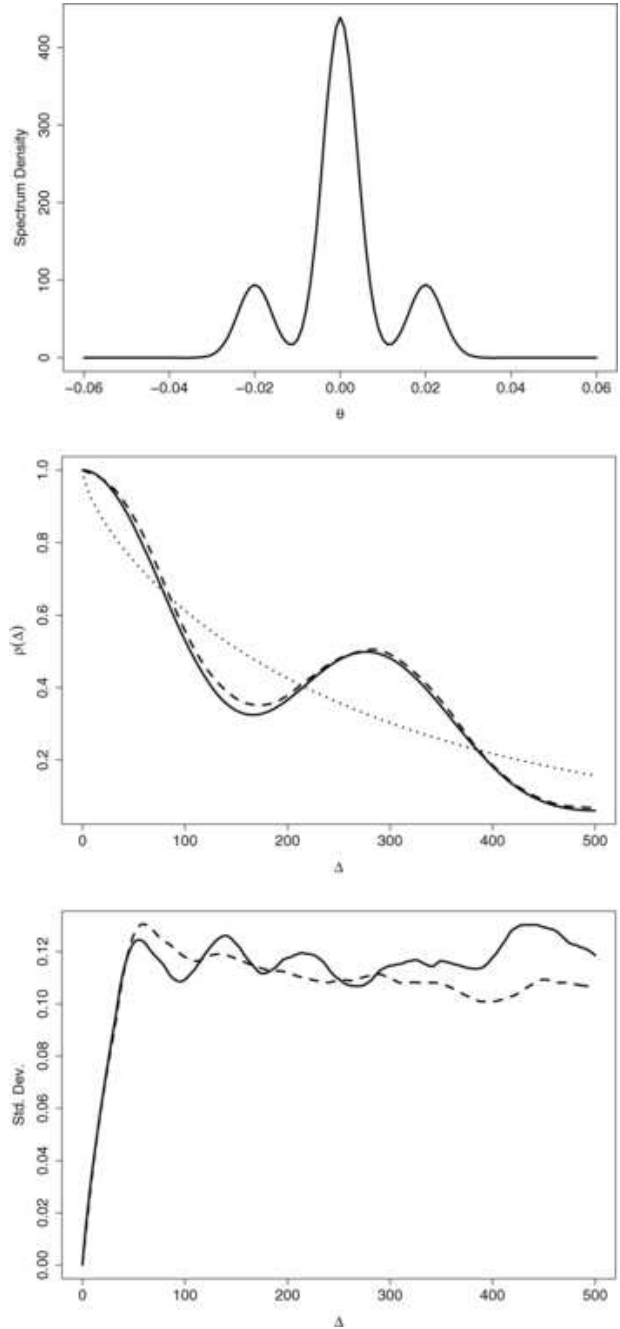

FIG. 5. *Simulation 2. Upper panel: the spectral density of the correlation used in the simulation; middle panel: the solid curve is the true correlation function, the dashed curve is the mean of the kernel correlation estimator and the dotted curve is the best Matérn approximation to the true correlation; lower panel: the solid curve is the true pointwise standard deviation for the kernel correlation estimator, the dashed curve is the mean for the bootstrap standard deviation estimator.*

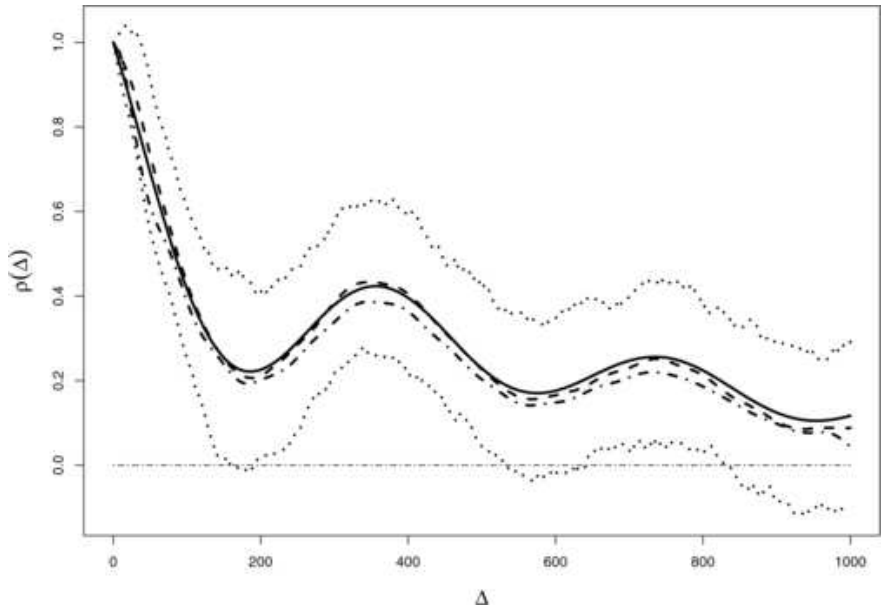

Fig. 6. *Simulation 3. The solid curve is the true correlation function, the dashed curve is the mean of our kernel estimator $\widehat{\rho}$, the two dotted curves are the 5% and 95% percentiles of $\widehat{\rho}$, and the dot–dash curve is the mean of the adjusted estimator $\widetilde{\rho}$.*

**6. Discussion.** We have proposed an estimator of stationary correlation functions for longitudinal or spatial data where within-subject observations have a complex data structure. The application we presented has a functional data flavor, in that each unit (crypt) in a "time series" has subunits (cells), the values from which can be viewed as a function. However, in this paper, we have focused on estimating the spatial correlation between the units.

We established an asymptotic normal limit distribution for the proposed estimator. The techniques used in our theoretical derivation were significantly different from those of the standard kernel regression literature. In our theoretical framework, as long as we have an increasing number of observations within a subject, each subject yields a consistent estimate of the correlation function. Our method and theory are especially useful to the cases where the number of subjects is limited but we have a relatively large number of repeated measurements within each subject. Since having more subjects will just further reduce the variation of the estimator, our main theorems hold when $R$ goes to infinity as well. In that case, we need to replace the condition that $Lh^5 = O(1)$ in assumption 7 in Section 3 with $RLh^5 = O(1)$. In fact, when the number of subjects $R \to \infty$, we can consistently estimate the within-subject covariance without a large number of units within each subject. For example, Yao, Müller and Wang [18] proposed using smoothing methods to estimate within-subject covariance for sparse longitudinal data.

In spatial statistics, many authors have considered the setup under the intrinsic stationary assumption (Besag, York and Mollié [2], Besag and Higdon [1]). This is weaker than our second-order stationary assumption. In our case, each unit within a subject has further structure, so that we can define a cross-variogram (Cressie [3]) instead of the covariance function $\mathcal{V}(x_1, x_2, \Delta)$, and similar limiting distribution theorems can be proved as in Theorem 1 and 2. However, when it comes to spatio-temporal modeling, many authors (Cressie and Huang [4], Stein [17]) would still focus their attention on covariance estimation because it is a more natural way to introduce the separable structure (3). In our data analysis, we provided some practical ideas to justify the separable structure in our data, where we compare the cross-validation scores with and without the separable assumption.

We proposed a weighted bootstrap method to estimate the standard deviation of the correlation estimator $\widehat{\rho}$, where the weights were constructed based on the outcome from Lemma A.2 in the proofs. Our simulation studies show that the proposed correlation estimator and the weighted bootstrap standard deviation estimator work well numerically for finite sample sizes.

## APPENDIX: PROOFS

The proofs are organized in the following way: in Section A.1, we provide lemmas regarding asymptotic properties of the covariance estimators when there is only one subject; in Section A.2, we provide lemmas on the estimators with multiple subjects, and the proofs of Theorems 1, 2 and Corollary 1 are given at the end.

**A.1. Estimation within one subject.** We first discuss the case where there is only one subject and the number of units goes to infinity. Let $N(\cdot)$ be the inhomogeneous Poisson process on $[0, L]$ with local intensity $\nu g^*(s)$. As in Karr [10], denote $N_2(ds_1, ds_2) = N(ds_1)N(ds_2)I(s_1 \neq s_2)$. Let $\Theta(s, \cdot)$ denote the unit-level mean at unit location $s$, and $\Psi(\cdot)$ denote the subject-level mean. Define

$$a(\Delta) = L^{-1} \sum_i \sum_{k \neq i} K_h(\Delta - \Delta_{ik})$$

$$= L^{-1} \int_0^L \int_0^L K_h\{\Delta - (s_1 - s_2)\} N_2(ds_1, ds_2);$$

$$b(x_1, x_2, \Delta) = L^{-1} \sum_i \sum_{k \neq i} K_h(\Delta - \Delta_{ik})\{Y(S_i, x_1) - \Psi(x_1)\}\{Y(S_k, x_2) - \Psi(x_2)\}$$

$$= L^{-1} \int_0^L \int_0^L K_h\{\Delta - (s_1 - s_2)\}\{Y(s_1, x_1) - \Psi(x_1)\}$$

$$\times \{Y(s_2, x_2) - \Psi(x_2)\} N_2(ds_1, ds_2).$$

LEMMA A.1. *Let $X_1$ and $X_2$ be real-valued random variables measurable with respect to $\mathcal{F}\{[0,t]\}$ and $\mathcal{F}\{[t+\tau,\infty)\}$ respectively, such that $|X_i| < C_i$, $i=1,2$. Then $|\text{cov}(X_1, X_2)| \leq 4C_1 C_2 \alpha(\tau)$. If $X_1$ and $X_2$ are complex random variables, this inequality holds with the constant 4 replaced by 16.*

PROOF. The proof is analogous to that of Theorem 17.2.1 in Ibragimov and Linnik [9]. Denote $T_1 = [0,t]$, $T_2 = [t+\tau, \infty)$. Then we have

$$\begin{aligned}
|E(X_1 X_2) &- E(X_1) E(X_2)| \\
&= |E[E\{X_1 X_2 | \mathcal{F}(T_1)\}] - E(X_1) E(X_2)| \\
&= |E(X_1 [E\{X_2 | \mathcal{F}(T_1)\} - E(X_2)])| \\
&\leq C_1 E |E\{X_2 | \mathcal{F}(T_1)\} - E(X_2)| \\
&= C_1 E(u_1 [E\{X_2 | \mathcal{F}(T_1)\} - E(X_2)]),
\end{aligned}$$

where $u_1 = \text{sign}[E\{X_2|\mathcal{F}(T_1)\} - E(X_2)]$. It is easy to see that $u_1$ is measurable with respect to $\mathcal{F}(T_1)$, and therefore $|E(X_1 X_2) - E(X_1) E(X_2)| \leq C_1 |E(u_1 X_2) - E(u_1) E(X_2)|$. By the same argument, we have $|E(u_1 X_2) - E(u_1) E(X_2)| \leq C_2 |E(u_1 u_2) - E(u_1) E(u_2)|$, where $u_2 = \text{sign}[E\{u_1|\mathcal{F}(T_2)\} - E(u_1)]$. Now we have $|E(X_1 X_2) - E(X_1) E(X_2)| \leq C_1 C_2 |E(u_1 u_2) - E(u_1) E(u_2)|$. Define the events $A_1 = \{u_1 = 1\} \in \mathcal{F}(T_1)$, $\overline{A}_1 = \{u_1 = -1\} \in \mathcal{F}(T_1)$, $A_2 = \{u_2 = 1\} \in \mathcal{F}(T_2)$ and $\overline{A}_2 = \{u_2 = -1\} \in \mathcal{F}(T_2)$. Then

$$\begin{aligned}
|E(u_1 u_2) &- E(u_1) E(u_2)| \\
&= |P(A_1 A_2) - P(A_1 \overline{A}_2) - P(\overline{A}_1 A_2) + P(\overline{A}_1 \overline{A}_2) \\
&\quad - P(A_1) P(A_2) + P(A_1) P(\overline{A}_2) + P(\overline{A}_1) P(A_2) - P(\overline{A}_1) P(\overline{A}_2)| \\
&\leq |P(A_1 A_2) - P(A_1) P(A_2)| + |P(A_1 \overline{A}_2) - P(A_1) P(\overline{A}_2)| \\
&\quad + |P(\overline{A}_1 A_2) - P(\overline{A}_1) P(A_2)| + |P(\overline{A}_1 \overline{A}_2) - P(\overline{A}_1) P(\overline{A}_2)| \\
&\leq 4\alpha(\tau).
\end{aligned}$$

Thus, the proof is completed for the real random variable case. If $X_1$ and $X_2$ are complex, we can apply the same arguments to the real and imaginary parts separately. □

LEMMA A.2. *With the assumptions stated in Section 3, for any fixed $\Delta$, we have $a(\Delta) \to \nu^2 f_1(0)$ in the $L^2$ sense, as $L \to \infty$.*

PROOF. Recall that by definition of $f_1(\cdot)$, if $X_1$ and $X_2$ are independent and identically distributed with density $g(\cdot)$, then $f_1(u) = \int g(t+u) g(t) \, dt$ is the density of $X_1 - X_2$. Thus, for fixed $\Delta$,

$$\begin{aligned}E\{a(\Delta)\} &= \nu^2 L^{-1} \int\int_{s_1 \neq s_2} K_h\{\Delta - (s_1 - s_2)\}g(s_1/L)g(s_2/L)\,ds_1\,ds_2 \\ &= \nu^2 L \int_0^1 \int_0^1 K_h\{\Delta - L(t_1 - t_2)\}g(t_1)g(t_2)\,dt_1\,dt_2 \\ &= \nu^2 L \int\int K_h(\Delta - Lu)g(t_2 + u)g(t_2)\,du\,dt_2 \\ &= \nu^2 L \int K_h(\Delta - Lu)f_1(u)\,du = \nu^2 \int K(v)f_1\{(\Delta - hv)/L\}\,dv \\ &= \nu^2 \int K(v)\{f_1(0) + O(L^{-1})\}\,dv = \nu^2 f_1(0) + O(L^{-1}).\end{aligned}$$

Next,

$$E\{a^2(\Delta)\} = L^{-2} \int_0^L \int_0^L \int_0^L \int_0^L K_h\{\Delta - (s_1 - s_2)\}K_h\{\Delta - (s_3 - s_4)\} \\ \times E\{N_2(ds_1, ds_2)N_2(ds_3, ds_4)\}.$$

Calculations as in Guan, Sherman and Calvin [6] show that

$$\begin{aligned}&E\{N_2(ds_1, ds_2)N_2(ds_3, ds_4)\} \\ &\quad = \nu^4 g^*(s_1)g^*(s_2)g^*(s_3)g^*(s_4)\,ds_1\,ds_2\,ds_3\,ds_4 \\ &\qquad + \nu^3 g^*(s_1)g^*(s_2)g^*(s_4)\varepsilon_{s_1}(ds_3)\,ds_1\,ds_2\,ds_4 \\ &\qquad + \nu^3 g^*(s_1)g^*(s_2)g^*(s_3)\varepsilon_{s_1}(ds_4)\,ds_1\,ds_2\,ds_3 \\ &\qquad + \nu^3 g^*(s_1)g^*(s_2)g^*(s_4)\varepsilon_{s_2}(ds_3)\,ds_1\,ds_2\,ds_4 \\ &\qquad + \nu^3 g^*(s_1)g^*(s_2)g^*(s_3)\varepsilon_{s_2}(ds_4)\,ds_1\,ds_2\,ds_3 \\ &\qquad + \nu^2 g^*(s_1)g^*(s_2)\varepsilon_{s_1}(ds_3)\varepsilon_{s_2}(ds_4)\,ds_1\,ds_2 \\ &\qquad + \nu^2 g^*(s_1)g^*(s_2)\varepsilon_{s_1}(ds_4)\varepsilon_{s_2}(ds_3)\,ds_1\,ds_2,\end{aligned}$$

where $\varepsilon_x(\cdot)$ is a point measure defined in Karr [10], such that $\varepsilon_x(dy) = 1$ if $x \in dy$, 0 otherwise. Here $dy$ is defined to be a small disc centered at $y$. There are seven terms in the expression above, so the expression for $E\{a^2(\Delta)\}$ can be decomposed into seven integrals; denote them as $A_{11}$–$A_{17}$. Similar to the calculations of $E\{a(\Delta)\}$, we have

$$\begin{aligned}A_{11} &= \nu^4 L^{-2} \int\int_{s_1 \neq s_2} \int\int_{s_3 \neq s_4} K_h\{\Delta - (s_1 - s_2)\}K_h\{\Delta - (s_3 - s_4)\} \\ &\qquad \times g(s_1/L)g(s_2/L)g(s_3/L)g(s_4/L)\,ds_1\,ds_2\,ds_3\,ds_4 \\ &= \nu^4 f_1^2(0) + o(1),\end{aligned}$$

$$A_{12} = \nu^3 L^{-2} \int_{s_1 \neq s_2, s_4} K_h\{\Delta - (s_1 - s_2)\} K_h\{\Delta - (s_1 - s_4)\}$$

$$\times g(s_1/L) g(s_2/L) g(s_4/L) \, ds_1 \, ds_2 \, ds_4$$

$$= \nu^3 L \int \int K_h(\Delta - Lu_1) K_h\{\Delta - L(u_1 - u_2)\} f_2(u_1, u_2) \, du_1 \, du_2$$

(by definition of $f_2$)

$$= \nu^3 L^{-1} \int \int K(v_1) K(v_2) f_2\{(\Delta - hv_1)/L, (v_2 - v_1)h/L\} \, dv_1 \, dv_2$$

$$= \nu^3 L^{-1} f_2(0, 0) + O(L^{-2}).$$

Similarly, $A_{13}$–$A_{15}$ are of order $O(L^{-1})$. Next,

$$A_{16} = \nu^2 L^{-2} \int_{s_1 \neq s_2} K_h^2\{\Delta - (s_1 - s_2)\} g(s_1/L) g(s_2/L) \, ds_1 \, ds_2$$

$$= \nu^2 \int K_h^2(\Delta - Lu) f_1(u) \, du$$

$$= \nu^2 L^{-1} h^{-1} \int K^2(v) f_1\{(\Delta - hv)/L\} \, dv$$

$$= \nu^2 L^{-1} h^{-1} f_1(0) R_K + o(Lh^{-1}).$$

Similarly, we can show that $A_{17}$ is of the same order as $A_{16}$. This means that $A_{11}$ is the leading term in $E\{a^2(\Delta)\}$. Hence, $E\{a(\Delta) - \nu^2 f_1(0)\}^2 \to 0$, completing the proof. $\square$

LEMMA A.3. *For any fixed $\Delta$, define $\beta(x_1, x_2, \Delta) = b(x_1, x_2, \Delta) - a(\Delta) \times \mathcal{V}(x_1, x_2, \Delta)$. Then*

$$E\{\beta(x_1, x_2, \Delta)\} = \nu^2 f_1(0)\{\mathcal{V}^{(0,0,2)}(x_1, x_2, \Delta)\sigma_K^2 h^2/2 + o(h^2)\},$$

$$\operatorname{cov}\{\beta(x_1, x_2, \Delta), \beta(x_3, x_4, \Delta')\}$$

$$= \nu^2 L^{-1} h^{-1} R_K f_1(0)$$

$$\times [I(\Delta = \Delta')\{\mathcal{M}(x_1, x_2, x_3, x_4, \Delta, \Delta, 0)$$

$$+ I(x_2 = x_4)\sigma_\varepsilon^2 \mathcal{V}(x_1, x_3, 0) + I(x_1 = x_3)\sigma_\varepsilon^2 \mathcal{V}(x_2, x_4, 0)$$

$$+ I(x_1 = x_3, x_2 = x_4)\sigma_\varepsilon^4\}$$

$$+ I(\Delta = -\Delta')\{\mathcal{M}(x_1, x_2, x_3, x_4, \Delta, -\Delta, -\Delta)$$

$$+ I(x_2 = x_3)\sigma_\varepsilon^2 \mathcal{V}(x_1, x_4, 0)$$

$$+ I(x_1 = x_4)\sigma_\varepsilon^2 \mathcal{V}(x_2, x_3, 0)$$

$$+ I(x_1 = x_4, x_2 = x_3)\sigma_\varepsilon^4\}]$$

$$+ o(L^{-1} h^{-1}),$$

where $\mathcal{V}^{(0,0,2)}(x_1, x_2, \Delta) = \partial^2 \mathcal{V}(x_1, x_2, \Delta)/\partial \Delta^2$.

PROOF. Rewrite
$$\beta(x_1, x_2, \Delta)$$
$$= L^{-1} \int\int K_h\{\Delta - (s_1 - s_2)\}$$
$$\times [\{Y(s_1, x_1) - \Psi(x_1)\}$$
$$\times \{Y(s_2, x_2) - \Psi(x_2)\} - \mathcal{V}(x_1, x_2, \Delta)]N_2(ds_1, ds_2).$$

It follows that
$$E\{\beta(x_1, x_2, \Delta)\}$$
$$= \nu^2 L^{-1} \int\int_{s_1 \neq s_2} K_h\{\Delta - (s_1 - s_2)\}$$
$$\times \{\mathcal{V}(x_1, x_2, s_1 - s_2) - \mathcal{V}(x_1, x_2, \Delta)\}$$
$$\times g(s_1/L)g(s_2/L)\, ds_1\, ds_2$$
$$= \nu^2 L \int K_h(\Delta - Lu)\{\mathcal{V}(x_1, x_2, Lu) - \mathcal{V}(x_1, x_2, \Delta)\}f_1(u)\, du$$
$$= \nu^2 \int K(v)\{-\mathcal{V}^{(0,0,1)}(x_1, x_2, \Delta)hv$$
$$+ \mathcal{V}^{(0,0,2)}(x_1, x_2, \Delta)h^2v^2/2 + o(h^2)\}\{f_1(0) + O(L^{-1})\}\, dv$$
$$= \nu^2\{f_1(0)\mathcal{V}^{(0,0,2)}(x_1, x_2, \Delta)\sigma_K^2 h^2/2 + o(h^2)\}.$$

In addition,
$$\text{cov}\{\beta(x_1, x_2, \Delta), \beta(x_3, x_4, \Delta')\}$$
$$= L^{-2} \int\int\int\int K_h\{\Delta - (s_1 - s_2)\}K_h\{\Delta' - (s_3 - s_4)\}$$
$$\times [\mathcal{V}(x_1, x_2, \Delta)\mathcal{V}(x_3, x_4, \Delta')$$
$$- \mathcal{V}(x_1, x_2, s_1 - s_2)\mathcal{V}(x_3, x_4, \Delta')$$
$$- \mathcal{V}(x_1, x_2, \Delta)\mathcal{V}(x_3, x_4, s_3 - s_4)$$
$$+ \mathcal{M}\{x_1, x_2, x_3, x_4, (s_1 - s_2), (s_3 - s_4), (s_2 - s_4)\}$$
$$+ \mathcal{V}(x_1, x_2, s_1 - s_2)\mathcal{V}(x_3, x_4, s_3 - s_4)$$
$$+ I(s_1 = s_3)I(s_2 \neq s_4)I(x_1 = x_3)$$
$$\times \sigma_\varepsilon^2 \mathcal{V}\{x_2, x_4, (s_2 - s_4)\}$$
$$+ I(s_1 = s_4)I(s_2 \neq s_3)I(x_1 = x_4)$$

$$\times \sigma_\varepsilon^2 \mathcal{V}\{x_2, x_3, (s_2 - s_3)\}$$
$$+ I(s_2 = s_3) I(s_1 \neq s_4) I(x_2 = x_3)$$
$$\times \sigma_\varepsilon^2 \mathcal{V}\{x_1, x_4, (s_1 - s_4)\}$$
$$+ I(s_2 = s_4) I(s_1 \neq s_3) I(x_2 = x_4)$$
$$\times \sigma_\varepsilon^2 \mathcal{V}\{x_1, x_3, (s_1 - s_3)\}$$
$$+ I(s_1 = s_3, s_2 = s_4)$$
$$\times \{I(x_2 = x_4) \sigma_\varepsilon^2 \mathcal{V}(x_1, x_3, 0)$$
$$+ I(x_1 = x_3) \sigma_\varepsilon^2 \mathcal{V}(x_2, x_4, 0)$$
$$+ I(x_1 = x_3, x_2 = x_4) \sigma_\varepsilon^4\}$$
$$+ I(s_1 = s_4, s_2 = s_3)$$
$$\times \{I(x_2 = x_3) \sigma_\varepsilon^2 \mathcal{V}(x_1, x_4, 0)$$
$$+ I(x_1 = x_4) \sigma_\varepsilon^2 \mathcal{V}(x_2, x_3, 0)$$
$$+ I(x_1 = x_4, x_2 = x_3) \sigma_\varepsilon^4\}]$$
$$\times E\{N_2(ds_1, ds_2) N_2(ds_3, ds_4)\}$$
$$- \nu^4 L^{-2} \int \int \int \int K_h\{\Delta - (s_1 - s_2)\} K_h\{\Delta' - (s_3 - s_4)\}$$
$$\times \{\mathcal{V}(x_1, x_2, s_1 - s_2) - \mathcal{V}(x_1, x_2, \Delta)\}$$
$$\times \{\mathcal{V}(x_3, x_4, s_3 - s_4) - \mathcal{V}(x_3, x_4, \Delta')\}$$
$$\times g(s_1/L) g(s_2/L) g(s_3/L) g(s_4/L) \, ds_1 \, ds_2 \, ds_3 \, ds_4.$$

As in Lemma A.2, according to the expression for $E\{N_2(ds_1, ds_2) N_2(ds_3, ds_4)\}$, we can summarize this covariance expression as the sum of seven terms, denoted as $A_{21}$–$A_{27}$. We have

$$A_{21} = \nu^4 L^{-2} \int_0^L \int_0^L \int_0^L \int_0^L K_h\{\Delta - (s_1 - s_2)\} K_h\{\Delta' - (s_3 - s_4)\}$$
$$\times \mathcal{M}\{x_1, x_2, x_1, x_2, (s_1 - s_2), (s_3 - s_4), (s_2 - s_4)\}$$
$$\times g(s_1/L) g(s_2/L) g(s_3/L) g(s_4/L) \, ds_1 \, ds_2 \, ds_3 \, ds_4$$
$$= \nu^4 L^2 \int_0^1 \int_0^1 \int_0^1 \int_0^1 K_h\{\Delta - L(t_1 - t_2)\} K_h\{\Delta' - L(t_3 - t_4)\}$$
$$\times \mathcal{M}\{x_1, x_2, x_1, x_2, L(t_1 - t_2), L(t_3 - t_4), L(t_2 - t_4)\}$$
$$\times g(t_1) g(t_2) g(t_3) g(t_4) \, dt_1 \, dt_2 \, dt_3 \, dt_4$$
$$= \nu^4 L^2 \int \int \int K_h(\Delta - L u_1) K_h(\Delta' - L u_2)$$

$$\times \mathcal{M}(x_1, x_2, x_1, x_2, Lu_1, Lu_2, Lu_3)$$
$$\times f_3(u_1, u_2, u_3) \, du_1 \, du_2 \, du_3$$
$$= \nu^4 L^{-1} \int \int \int K(v_1) K(v_2) \mathcal{M}(x_1, x_2, x_1, x_2, \Delta - hv_1, \Delta' - hv_2, v_3)$$
$$\times f_3\{(\Delta - hv_1)/L, (\Delta' - hv_2)/L, v_3/L\} \, dv_1 \, dv_2 \, dv_3$$
$$\leq \nu^4 L^{-1} C \int \mathcal{M}(x_1, x_2, x_1, x_2, \Delta, \Delta', v) \, dv + o(L^{-1}),$$

where $C$ is the upper bound for the density function $f_3(u, v, w)$ on $[-1, 1]^3$. By assumption 1 in Section 3 that $g(\cdot)$ is bounded, one can easily derive that $C$ is a finite constant. The second term is

$$A_{22} = \nu^3 L^{-2} \int \int \int K_h\{\Delta - (s_1 - s_2)\} K_h\{\Delta' - (s_1 - s_4)\}$$
$$\times ([\mathcal{V}(x_1, x_2, \Delta) - \mathcal{V}\{x_1, x_2, (s_1 - s_2)\}]$$
$$\times [\mathcal{V}(x_3, x_4, \Delta') - \mathcal{V}\{x_3, x_4, (s_1 - s_4)\}]$$
$$+ \mathcal{M}\{x_1, x_2, x_3, x_4, (s_1 - s_2), (s_1 - s_4), (s_2 - s_4)\}$$
$$+ I(x_1 = x_3) \sigma_\varepsilon^2 \mathcal{V}\{x_2, x_4, (s_2 - s_4)\})$$
$$\times g(s_1/L) g(s_2/L) g(s_4/L) \, ds_1 \, ds_2 \, ds_4$$
$$= \nu^3 L \int \int \int K_h\{\Delta - L(t_1 - t_2)\} K_h\{\Delta' - L(t_1 - t_4)\}$$
$$\times ([\mathcal{V}(x_1, x_2, \Delta) - \mathcal{V}\{x_1, x_2, L(t_1 - t_2)\}]$$
$$\times [\mathcal{V}(x_3, x_4, \Delta') - \mathcal{V}\{x_3, x_4, L(t_1 - t_4)\}]$$
$$+ \mathcal{M}\{x_1, x_2, x_3, x_4, L(t_1 - t_2), L(t_1 - t_4), L(t_2 - t_4)\}$$
$$+ I(x_1 = x_3) \sigma_\varepsilon^2 \mathcal{V}\{x_2, x_4, L(t_2 - t_4)\})$$
$$\times g(t_1) g(t_2) g(t_4) \, dt_1 \, dt_2 \, dt_4$$
$$= \nu^3 L \int \int K_h(\Delta + Lu_1) K_h(\Delta' + Lu_2)$$
$$\times [\{\mathcal{V}(x_1, x_2, \Delta) - \mathcal{V}(x_1, x_2, -Lu_1)\}$$
$$\times \{\mathcal{V}(x_3, x_4, \Delta') - \mathcal{V}(x_3, x_4, -Lu_2)\}$$
$$+ \mathcal{M}\{x_1, x_2, x_3, x_4, -Lu_1, -Lu_2, L(u_1 - u_2)\}$$
$$+ I(x_1 = x_3) \sigma_\varepsilon^2 \mathcal{V}\{x_2, x_4, L(u_1 - u_2)\}]$$
$$\times f_2(u_1, u_2) \, du_1 \, du_2$$
$$= \nu^3 L^{-1} \int \int K(v_1) K(v_2) [I(x_1 = x_3) \sigma_\varepsilon^2 \mathcal{V}\{x_2, x_4, (v_1 - v_2)h + \Delta' - \Delta\}$$

$$+ \{\mathcal{V}(x_1, x_2, \Delta) - \mathcal{V}(x_1, x_2, \Delta - hv_1)\}$$
$$+ \{\mathcal{V}(x_3, x_4, \Delta) - \mathcal{V}(x_3, x_4, \Delta' - hv_2)\}$$
$$+ \mathcal{M}\{x_1, x_2, x_3, x_4, \Delta - hv_1, \Delta' - hv_2,$$
$$(v_1 - v_2)h + \Delta' - \Delta\}]$$
$$\times f_2\{(-\Delta + hv_1)/L, (-\Delta + v_1 h)/L\} \, dv_1 \, dv_2$$
$$= \nu^3 L^{-1} f_2(0,0) \{\mathcal{M}(x_1, x_2, x_3, x_4, \Delta, \Delta', \Delta' - \Delta)$$
$$+ I(x_1 = x_3)\sigma_\varepsilon^2 \mathcal{V}(x_2, x_4, \Delta' - \Delta)\} + o(L^{-1}).$$

It is easy to see that $A_{23}$–$A_{25}$ have the same order as $A_{22}$. Further, we have

$$A_{26} = \nu^2 L^{-2} \int\int K_h\{\Delta - (s_1 - s_2)\} K_h\{\Delta' - (s_1 - s_2)\}$$
$$\times (\mathcal{M}\{x_1, x_2, x_3, x_4, (s_1 - s_2), (s_1 - s_2), 0\}$$
$$+ [\mathcal{V}(x_1, x_2, \Delta) - \mathcal{V}\{x_1, x_2, (s_1 - s_2)\}]$$
$$\times [\mathcal{V}(x_3, x_4, \Delta') - \mathcal{V}\{x_3, x_4, (s_1 - s_2)\}]$$
$$+ \{I(x_2 = x_4)\sigma_\varepsilon^2 \mathcal{V}(x_1, x_3, 0)$$
$$+ I(x_1 = x_3)\sigma_\varepsilon^2 \mathcal{V}(x_2, x_4, 0) + I(x_1 = x_3, x_2 = x_4)\sigma_\varepsilon^4\})$$
$$\times g(s_1/L)g(s_2/L) \, ds_1 \, ds_2$$
$$= I(\Delta = \Delta')\nu^2 \int\int K_h^2\{\Delta - L(t_1 - t_2)\}$$
$$\times (\mathcal{M}\{x_1, x_2, x_3, x_4, L(t_1 - t_2), L(t_1 - t_2), 0\}$$
$$+ [\mathcal{V}(x_1, x_2, \Delta) - \mathcal{V}\{x_1, x_2, L(t_1 - t_2)\}]$$
$$\times [\mathcal{V}(x_3, x_4, \Delta) - \mathcal{V}\{x_3, x_4, L(t_1 - t_2)\}]$$
$$+ \{I(x_2 = x_4)\sigma_\varepsilon^2 \mathcal{V}(x_1, x_3, 0)$$
$$+ I(x_1 = x_3)\sigma_\varepsilon^2 \mathcal{V}(x_2, x_4, 0)$$
$$+ I(x_1 = x_3, x_2 = x_4)\sigma_\varepsilon^4\})$$
$$\times g(t_1)g(t_2) \, dt_1 \, dt_2$$
$$= I(\Delta = \Delta')\nu^2 \int K_h^2(\Delta - Lu)$$
$$\times [\mathcal{M}(x_1, x_2, x_3, x_4, Lu, Lu, 0)$$
$$+ \{\mathcal{V}(x_1, x_2, \Delta) - \mathcal{V}(x_1, x_2, Lu)\}$$
$$\times \{\mathcal{V}(x_3, x_4, \Delta) - \mathcal{V}(x_3, x_4, Lu)\}$$
$$+ \{I(x_2 = x_4)\sigma_\varepsilon^2 \mathcal{V}(x_1, x_3, 0)$$

$$+ I(x_1 = x_3)\sigma_\varepsilon^2 \mathcal{V}(x_2, x_4, 0)$$
$$+ I(x_1 = x_3, x_2 = x_4)\sigma_\varepsilon^4\}]$$
$$\times f_1(u)\,du$$
$$= I(\Delta = \Delta')\nu^2 L^{-1} h^{-1} \int K^2(v)[\mathcal{M}(x_1, x_2, x_3, x_4, \Delta - hv, \Delta - hv, 0)$$
$$+ \{\mathcal{V}(x_1, x_2, \Delta) - \mathcal{V}(x_1, x_2, \Delta - hv)\}$$
$$\times \{\mathcal{V}(x_3, x_4, \Delta) - \mathcal{V}(x_3, x_4, \Delta - hv)\}$$
$$+ \{I(x_2 = x_4)\sigma_\varepsilon^2 \mathcal{V}(x_1, x_3, 0)$$
$$+ I(x_1 = x_3)\sigma_\varepsilon^2 \mathcal{V}(x_2, x_4, 0)$$
$$+ I(x_1 = x_3, x_2 = x_4)\sigma_\varepsilon^4\}]$$
$$\times f_1\{(\Delta - hv)/L\}\,dv$$
$$= I(\Delta = \Delta')\nu^2 L^{-1} h^{-1} R_K f_1(0)$$
$$\times [\mathcal{M}(x_1, x_2, x_3, x_4, \Delta, \Delta, 0)$$
$$+ \{I(x_2 = x_4)\sigma_\varepsilon^2 \mathcal{V}(x_1, x_3, 0) + I(x_1 = x_3)\sigma_\varepsilon^2 \mathcal{V}(x_2, x_4, 0)$$
$$+ I(x_1 = x_3, x_2 = x_4)\sigma_\varepsilon^4\} + o(1)].$$

Similarly,
$$A_{27} = \nu^2 L^{-2} \int\int K_h\{\Delta - (s_1 - s_2)\} K_h\{\Delta' - (s_2 - s_1)\}$$
$$\times (\mathcal{M}\{x_1, x_2, x_3, x_4, (s_1 - s_2), (s_2 - s_1), (s_2 - s_1)\}$$
$$+ [\mathcal{V}(x_1, x_2, \Delta) - \mathcal{V}\{x_1, x_2, (s_1 - s_2)\}]$$
$$\times [\mathcal{V}(x_3, x_4, \Delta') - \mathcal{V}\{x_3, x_4, (s_2 - s_1)\}]$$
$$+ \{I(x_2 = x_3)\sigma_\varepsilon^2 \mathcal{V}(x_1, x_4, 0) + I(x_1 = x_4)\sigma_\varepsilon^2 \mathcal{V}(x_2, x_3, 0)$$
$$+ I(x_1 = x_4, x_2 = x_3)\sigma_\varepsilon^4\})$$
$$\times g(s_1/L)g(s_2/L)\,ds_1\,ds_2$$
$$= I(\Delta = -\Delta')\nu^2 \int\int K_h^2\{\Delta - L(t_1 - t_2)\}$$
$$\times (\mathcal{M}\{x_1, x_2, x_3, x_4, L(t_1 - t_2), L(t_2 - t_1), L(t_2 - t_1)\}$$
$$+ [\mathcal{V}(x_1, x_2, \Delta) - \mathcal{V}\{x_1, x_2, L(t_1 - t_2)\}]$$
$$\times [\mathcal{V}(x_3, x_4, -\Delta) - \mathcal{V}\{x_3, x_4, L(t_2 - t_1)\}]$$
$$+ \{I(x_2 = x_3)\sigma_\varepsilon^2 \mathcal{V}(x_1, x_4, 0)$$
$$+ I(x_1 = x_4)\sigma_\varepsilon^2 \mathcal{V}(x_2, x_3, 0)$$

$$+ I(x_1 = x_4, x_2 = x_3)\sigma_\varepsilon^4\})$$
$$\times g(t_1)g(t_2)\,dt_1\,dt_2$$
$$= I(\Delta = -\Delta')\nu^2 \int K_h^2(\Delta - Lu)[\mathcal{M}(x_1, x_2, x_3, x_4, Lu, -Lu, -Lu)$$
$$+ \{\mathcal{V}(x_1, x_2, \Delta) - \mathcal{V}(x_1, x_2, Lu)\}$$
$$\times \{\mathcal{V}(x_3, x_4, \Delta) - \mathcal{V}(x_3, x_4, Lu)\}$$
$$+ \{I(x_2 = x_3)\sigma_\varepsilon^2 \mathcal{V}(x_1, x_4, 0)$$
$$+ I(x_1 = x_4)\sigma_\varepsilon^2 \mathcal{V}(x_2, x_3, 0)$$
$$+ I(x_1 = x_4, x_2 = x_3)\sigma_\varepsilon^4\}]$$
$$\times f_1(u)\,du$$
$$= I(\Delta = -\Delta')\nu^2 L^{-1} h^{-1} \int K^2(v)$$
$$\times [\mathcal{M}(x_1, x_2, x_3, x_4, \Delta - hv, -\Delta + hv, -\Delta + hv)$$
$$+ \{\mathcal{V}(x_1, x_2, \Delta) - \mathcal{V}(x_1, x_2, \Delta - hv)\}$$
$$\times \{\mathcal{V}(x_3, x_4, \Delta) - \mathcal{V}(x_3, x_4, \Delta - hv)\}$$
$$+ \{I(x_2 = x_3)\sigma_\varepsilon^2 \mathcal{V}(x_1, x_4, 0)$$
$$+ I(x_1 = x_4)\sigma_\varepsilon^2 \mathcal{V}(x_2, x_3, 0)$$
$$+ I(x_1 = x_4, x_2 = x_3)\sigma_\varepsilon^4\}]$$
$$\times f_1\{(\Delta - hv)/L\}\,dv$$
$$= I(\Delta = -\Delta')\nu^2 L^{-1} h^{-1} R_K f_1(0)$$
$$\times [\mathcal{M}(x_1, x_2, x_3, x_4, \Delta, -\Delta, -\Delta)$$
$$+ \{I(x_2 = x_3)\sigma_\varepsilon^2 \mathcal{V}(x_1, x_4, 0)$$
$$+ I(x_1 = x_4)\sigma_\varepsilon^2 \mathcal{V}(x_2, x_3, 0)$$
$$+ I(x_1 = x_4, x_2 = x_3)\sigma_\varepsilon^4\} + o(1)].$$

Both $A_{26}$ and $A_{27}$ are of order $O\{(Lh)^{-1}\}$, while the rest of the terms are of order $O(L^{-1})$. The proof is completed by summarizing the contribution of each term to $\mathrm{cov}\{\beta(x_1, x_2, \Delta), \beta(x_3, x_4, \Delta')\}$. □

LEMMA A.4. *With $\beta(x_1, x_2, \Delta)$ defined as in Lemma A.3, and with all the assumptions in Section 3, we have*

$$(Lh)^{1/2}[\beta(x_1, x_2, \Delta) - E\{\beta(x_1, x_2, \Delta)\}] \Rightarrow \mathrm{Normal}\{0, \nu^2 f_1(0)\sigma^2(x_1, x_2, \Delta)\},$$

where $\sigma^2(x_1,x_2,\Delta) = R_K\{\mathcal{M}(x_1,x_2,x_1,x_2,\Delta,\Delta,0)+\sigma_\varepsilon^2\mathcal{V}(x_1,x_1,0)+\sigma_\varepsilon^2\mathcal{V}(x_2,x_2,0) + \sigma_\varepsilon^4\} + I(\Delta=0)R_K[\{\mathcal{M}(x_1,x_2,x_1,x_2,0,0,0) + I(x_1=x_2)\{2\sigma_\varepsilon^2\mathcal{V}(x_1,x_1,0) + \sigma_\varepsilon^4\}]$.

PROOF. The proof has similar structure to that of Theorem 2 in Guan, Sherman and Calvin [6]. Define $a_1 = 0$, $b_1 = L^p - L^q$, $a_i = a_{i-1} + L^p$, $b_i = a_i + L^p - L^q$, $i = 2,\ldots,k_L$, for some $1/(1+\delta) < q < p < 1$ [$\delta$ is defined in (10)]. We thus have divided the interval $[0,L]$ into $k_L \approx L/L^p$ disjoint subintervals each having length $L^p - L^q$ and at least $L^q$ apart. Define $I_i = [a_i,b_i]$, $I = \bigcup_{i=1}^{k_L} I_i$, $I_i' = [a_i/L, b_i/L]$, $I' = \bigcup_{i=1}^{k_L} I_i'$ and

$$\beta_i(x_1,x_2,\Delta) = L^{-1}\int\int_{I_i\times I_i} K_h\{\Delta - (s_1-s_2)\}$$
$$\times [\{Y(s_1,x_1) - \Psi(x_1)\}\{Y(s_2,x_2) - \Psi(x_2)\}$$
$$- \mathcal{V}(x_1,x_2,\Delta)]$$
$$\times N_2(ds_1,ds_2),$$

$$\widetilde{\beta}(x_1,x_2,\Delta) = \sum_{i=1}^{k_L}\beta_i(x_1,x_2,\Delta).$$

Define independent random variables $\gamma_i(x_1,x_2,\Delta)$ on a different probability space such that they have the same distributions as $\beta_i(x_1,x_2,\Delta)$, and define $\gamma(x_1,x_2,\Delta) = \sum_{i=1}^{k_L}\gamma_i(x_1,x_2,\Delta)$. Let $\phi(\xi)$ and $\psi(\xi)$ be the characteristic functions of $(Lh)^{1/2}[\widetilde{\beta}(x_1,x_2,\Delta) - E\{\widetilde{\beta}(x_1,x_2,\Delta)\}]$ and $(Lh)^{1/2}[\gamma(x_1,x_2,\Delta) - E\{\gamma(x_1,x_2,\Delta)\}]$, respectively.

We finish the proof in the following three steps:

(i) $([\beta(x_1,x_2,\Delta) - E\{\beta(x_1,x_2,\Delta)\}] - [\widetilde{\beta}(x_1,x_2,\Delta) - E\{\widetilde{\beta}(x_1,x_2,\Delta)\}]) \xrightarrow{p} 0$;
(ii) $\psi(\xi) - \phi(\xi) \to 0$;
(iii) $(Lh)^{1/2}[\gamma(x_1,x_2,\Delta) - E\{\gamma(x_1,x_2,\Delta)\}] \Rightarrow \text{Normal}\{0,\nu^2 f_1(0)\sigma^2(x_1,x_2,\Delta)\}$.

To show (i), notice that, with $|I_i| \to \infty$, calculations as in Lemma A.3 show that

$$\sum_{i=1}^{k_L}\text{var}\{\beta_i(x_1,x_2,\Delta)\} = \sum_{i=1}^{k_L}\nu^2 L^{-1}h^{-1}R_K f_{i,1}(0)\{\sigma^2(x_1,x_2,\Delta) + o(1)\},$$
(14)
where $f_{i,1}(u) = \int g_i(u+t)g_i(t)\,dt$ is the counterpart of $f_1(u)$, with $g_i(t) = g(t)I(t \in I_i')$. Since $g(\cdot)$ is bounded away from both 0 and $\infty$, $f_{i,1}(0) = \int_{I_i'} g^2(t) = O(|I_i'|) = O(L^{p-1})$ and $\text{var}\{\beta_i(x_1,x_2,\Delta)\} = O(L^{p-2}h^{-1})$.

SPATIAL CORRELATION FUNCTIONS 31Observe that $|I'| = \sum_{i=1}^{k_L} |I'_i| = k_L \times (L^p - L^q)/L \approx L/L^p \times (L^p - L^q)/L = 1 - L^{q-p} \to 1$, and

$$(15) \quad \sum_{i=1}^{k_L} f_{i,1}(0) = \sum_{i=1}^{k_L} \int_{I'_i} g(t)^2\, dt = \int_{I'} g(t)^2\, dt \to \int_0^1 g(t)^2\, dt = f_1(0).$$

Therefore, $\sum_{i=1}^{k_L} \mathrm{var}\{\beta_i(x_1, x_2, \Delta)\} = \mathrm{var}\{\beta(x_1, x_2, \Delta)\} + o(L^{-1}h^{-1})$. Further but equivalent derivations show that $\sum_{i \neq j} \mathrm{cov}\{\beta_i(x_1, x_2, \Delta), \beta_j(x_1, x_2, \Delta)\} = O(L^{-1})$. The calculations here are similar to those in Lemma A.3, except that the $i \neq j$ condition excluded terms like $A_{22}$ through $A_{27}$. Now we have

$$\mathrm{var}\{\widetilde{\beta}(x_1, x_2, \Delta)\} = \sum_{i=1}^{k_L} \mathrm{var}\{\beta_i(x_1, x_2, \Delta)\} + \sum_{i \neq j} \mathrm{cov}\{\beta_i(x_1, x_2, \Delta), \beta_j(x_1, x_2, \Delta)\}$$

$$= \mathrm{var}\{\beta(x_1, x_2, \Delta)\} + o(L^{-1}h^{-1}).$$

Similarly, one can show that

$$\mathrm{cov}\{\widetilde{\beta}(x_1, x_2, \Delta), \beta(x_1, x_2, \Delta)\} = \mathrm{var}\{\beta(x_1, x_2, \Delta)\} + o(L^{-1}h^{-1}).$$

Therefore, $(Lh)\,\mathrm{var}[\{\beta(x_1, x_2, \Delta) - \{\widetilde{\beta}(x_1, x_2, \Delta)\}] \to 0$, and step (i) is established. To show (ii), we follow similar arguments that prove Theorem 2 (S2) in Guan, Sherman and Calvin [6]. Denote $U_i = \exp(Ix(Lh)^{1/2}[\beta_i(x_1, x_2, \Delta) - E\{\beta_i(x_1, x_2, \Delta)\}])$, where $I$ is the unit imaginary number. Then by definition, $\phi(x) = E(\prod_{i=1}^{k_L} U_i)$, $\psi(x) = \prod_{i=1}^{k_L} E(U_i)$.

Observing $|E(U_i)| \leq 1$ for all $U_i$, we have

$|\phi(x) - \psi(x)|$

$$\leq \left| E\left(\prod_{i=1}^{k_L} U_i\right) - E\left(\prod_{i=1}^{k_L-1} U_i\right) E(U_{k_L}) \right| + \left| E\left(\prod_{i=1}^{k_L-1} U_i\right) E(U_{k_L}) - \prod_{i=1}^{k_L} E(U_i) \right|$$

$$\leq \left| E\left(\prod_{i=1}^{k_L} U_i\right) - E\left(\prod_{i=1}^{k_L-1} U_i\right) E(U_{k_L}) \right|$$

$$+ \left| E\left(\prod_{i=1}^{k_L-1} U_i\right) - \prod_{i=1}^{k_L-1} E(U_i) \right| |E(U_{k_L})|$$

$$\leq \left| E\left(\prod_{i=1}^{k_L} U_i\right) - E\left(\prod_{i=1}^{k_L-1} U_i\right) E(U_{k_L}) \right| + \left| E\left(\prod_{i=1}^{k_L-1} U_i\right) - \prod_{i=1}^{k_L-1} E(U_i) \right|.$$

By induction,

$$|\phi(x) - \psi(x)| \leq \sum_{j=1}^{k_L-1} \left| E\left(\prod_{i=1}^{j+1} U_i\right) - E\left(\prod_{i=1}^{j} U_i\right) E(U_{j+1}) \right|$$

$$= \sum_{j=1}^{k_L-1} \left| \text{cov}\left(\prod_{i=1}^{j} U_i, U_{j+1}\right) \right|.$$

Observe that $\prod_{i=1}^{j} U_i$ and $U_{j+1}$ are $\mathcal{F}([0, b_j])$- and $\mathcal{F}([a_{j+1}, b_{j+1}])$-measurable, respectively, with $|\prod_{i=1}^{j} U_i| \leq 1$ and $|U_{j+1}| \leq 1$, and the index sets are at least $L^q$ away. By Lemma A.1,

$$|\phi(x) - \psi(x)| \leq \sum_{j=1}^{k_L-1} 16\alpha(L^q) \leq 16L^{1-p} \times L^{-q\delta}.$$

By our choice of $p$ and $q$, it is easy to check $1 - p - q\delta < 0$, and therefore $|\phi(x) - \psi(x)| \to 0$.

(iii) can be proved by applying Lyapounov's central limit theorem and by the fact that

$$(Lh) \sum_{i=1}^{k_L} \text{var}\{\gamma_i(x_1, x_2, \Delta)\} \to \nu^2 f_1(0)\sigma^2(x_1, x_2, \Delta),$$

which has been shown in (14) and (15).

It remains to check Lyapounov's condition. By condition (9),

$$\sum_{i=1}^{k_L} \frac{E(|\gamma_i(x_1, x_2, \Delta) - E\{\gamma_i(x_1, x_2, \Delta)\}|^{2+\eta})}{[\text{var}\{\gamma(x_1, x_2, \Delta)\}]^{(2+\eta)/2}}$$

$$= L^{1-p} \frac{O\{(L^{p-2}h^{-1})^{(2+\eta)/2}\}}{O\{(L^{-1}h^{-1})^{(2+\eta)/2}\}}$$

$$= O(L^{-(1-p)\eta/2}) \to 0.$$

The proof is thus complete. □

LEMMA A.5. *Let $\vec{\beta}(\Delta)$ be the vector collecting all $\beta(x_1, x_2, \Delta)$ for distinct pairs of $(x_1, x_2)$. Then, with all assumptions above, for $\Delta' \neq \Delta$,*

$$(Lh)^{1/2} \begin{bmatrix} \vec{\beta}(\Delta) - E\{\vec{\beta}(\Delta)\} \\ \vec{\beta}(\Delta') - E\{\vec{\beta}(\Delta')\} \end{bmatrix}$$

$$\Rightarrow \text{Normal}\left\{0, \nu^2 f_1(0) \begin{pmatrix} \Sigma(\Delta) & C(\Delta, \Delta') \\ C^T(\Delta, \Delta') & \Sigma(\Delta') \end{pmatrix}\right\},$$

*where $\Sigma(\Delta)$ is the covariance matrix with the entry corresponding to $\text{cov}\{\beta(x_1, x_2, \Delta), \beta(x_3, x_4, \Delta)\}$ equal to $R_K\{\mathcal{M}(x_1, x_2, x_3, x_4, \Delta, \Delta, 0) + I(x_2 = x_4)\sigma_\varepsilon^2 \times \mathcal{V}(x_1, x_3, 0) + I(x_1 = x_3)\sigma_\varepsilon^2 \mathcal{V}(x_2, x_4, 0) + (x_1 = x_3, x_2 = x_4)\sigma_\varepsilon^4\} + I(\Delta = 0) \times R_K\{\mathcal{M}(x_1, x_2, x_3, x_4, 0, 0, 0) + I(x_1 = x_4)\sigma_\varepsilon^2 \mathcal{V}(x_2, x_3, 0) + I(x_2 = x_3)\sigma_\varepsilon^2 \mathcal{V}(x_1, x_4, 0) + I(x_1 = x_4, x_2 = x_3)\sigma_\varepsilon^4\}$ and $C(\Delta, \Delta')$ is the matrix with the entry corresponding to $\text{cov}\{\beta(x_1, x_2, \Delta), \beta(x_3, x_4, \Delta')\}$ equal to $I(\Delta' = -\Delta)\{\mathcal{M}(x_1, x_2,$* $x_3, x_4, \Delta, -\Delta, -\Delta) + I(x_2 = x_3)\sigma_\varepsilon^2 \mathcal{V}(x_1, x_4, 0) + I(x_1 = x_4)\sigma_\varepsilon^2 \mathcal{V}(x_2, x_3, 0) + I(x_1 = x_4, x_2 = x_3)\sigma_\varepsilon^4\}.$

PROOF. Using similar proofs as for Lemmas A.3 and A.4, we can show that any linear combination $\sum_{i=1}^{k} c_i \beta(x_{i1}, x_{i2}, \Delta) + \sum_{i=1}^{k'} c'_i \beta(x_{i1}, x_{i2}, \Delta')$ is asymptotically normal. By the Crámer–Wold device (Serfling [14]), the joint normality is established. □

NOTE. If $\Delta' = -\Delta$, the limiting distribution on the right-hand side is a degenerate multivariate normal distribution, because $\beta(x_1, x_1, \Delta) = \beta(x_1, x_1, -\Delta)$ for all $x_1$.

**A.2. Estimation with multiple subjects.** Now suppose we have $R$ subjects, and $R$ is a fixed number. Define

$$Y_{r,ik}(x_j, x_l) = \{Y_{rij} - \Psi_r(x_j)\}\{Y_{rkl} - \Psi_r(x_l)\},$$

$$a_r(\Delta) = L^{-1} \sum_i \sum_{k \neq i} K_h(\Delta - \Delta_{r,ik}),$$

$$b_r(x_j, x_l, \Delta) = L^{-1} \sum_i \sum_{k \neq i} Y_{r,ik}(x_j, x_l) K_h\{\Delta - \Delta_r(i,k)\},$$

$$\beta_r(x_j, x_l, \Delta) = b_r(x_j, x_l, \Delta) - a_r(\Delta)\mathcal{V}(x_j, x_l, \Delta),$$

$$c_r(x_j, \Delta) = L^{-1} \sum_i \sum_{k \neq i} \{Y_{rij} - \Psi_r(x_j)\} K_h\{\Delta - \Delta_r(i,k)\}.$$

Further, define $a(\Delta) = \sum_r a_r(\Delta)$, $b(x_j, x_l, \Delta) = \sum_r b_r(x_j, x_l, \Delta)$, $\beta(x_j, x_l, \Delta) = \sum_r \beta_r(x_j, x_l, \Delta)$ and $\widehat{\mathcal{V}}_0(x_1, x_2, \Delta) = b(x_1, x_2, \Delta)/a(\Delta)$. Let $\widehat{\mathcal{V}}_0(\Delta)$ and $\mathcal{V}(\Delta)$ be the vectors collecting all $\widehat{\mathcal{V}}_0(x_1, x_2, \Delta)$ and $\mathcal{V}(x_1, x_2, \Delta)$ for all distinct pairs of $(x_1, x_2)$, respectively.

LEMMA A.6. *With the assumptions in Section 3, for $\Delta' \neq \Delta$,*

$$(RLh)^{1/2} \left\{ \begin{array}{l} \widehat{\mathcal{V}}_0(\Delta) - \mathcal{V}(\Delta) - \sigma_K^2 \mathcal{V}^{(2)}(\Delta) h^2/2 \\ \widehat{\mathcal{V}}_0(\Delta') - \mathcal{V}(\Delta') - \sigma_K^2 \mathcal{V}^{(2)}(\Delta') h^2/2 \end{array} \right\}$$

$$\Rightarrow \text{Normal} \left[ 0, \{\nu^2 f_1(0)\}^{-1} \begin{pmatrix} \Sigma(\Delta) & C(\Delta, \Delta') \\ C^T(\Delta, \Delta') & \Sigma(\Delta') \end{pmatrix} \right],$$

*where $\mathcal{V}^{(2)}(\Delta)$ is the vector collecting $\mathcal{V}^{(0,0,2)}(x_1, x_2, \Delta)$ for all distinct pairs of $(x_1, x_2)$.*

PROOF. Notice that

$$\widehat{\mathcal{V}}_0(x_1, x_2, \Delta) - \mathcal{V}(x_1, x_2, \Delta)$$

$$= \left[\sum_{r=1}^{R}\{b_r(x_1, x_2, \Delta) - a_r(\Delta)\mathcal{V}(x_1, x_2, \Delta)\}\right] \Big/ \left\{\sum_{r=1}^{R} a_r(\Delta)\right\}$$
$$= \beta(x_1, x_2, \Delta)/a(\Delta).$$

Since subjects are independent, by Lemma A.2, $a(\Delta)/\{\nu^2 R f_1(0)\} \xrightarrow{p} 1$. Also, by Lemma A.5, $(R^{-1}Lh)^{1/2}\{\vec{\beta}(\Delta)^T, \vec{\beta}(\Delta')^T\}^T$ is asymptotically normal with covariance matrix given in Lemma A.5. Thus, by Slutsky's theorem [14],

$$(RLh)^{1/2}\begin{bmatrix} \beta(\Delta)/a(\Delta) - E\{\beta(\Delta)\}/a(\Delta) \\ \beta(\Delta')/a(\Delta') - E\{\beta(\Delta')\}/a(\Delta') \end{bmatrix}$$
$$\Rightarrow \text{Normal}\left[0, \{\nu^2 f_1(0)\}^{-1}\begin{pmatrix} \Sigma(\Delta) & C(\Delta, \Delta') \\ C^T(\Delta, \Delta') & \Sigma(\Delta') \end{pmatrix}\right].$$

Finally, by Lemma A.3, $E\{\beta(x_1, x_2, \Delta)\} = R\nu^2 f_1(0)\{\mathcal{V}^{(0,0,2)}(x_1, x_2, \Delta)\sigma_K^2 h^2/2 + o(h^2)\}$, so that we have $E\{\beta(x_1, x_2, \Delta)\}/a(\Delta) = \sigma_K^2 \mathcal{V}^{(0,0,2)}(x_1, x_2, \Delta)h^2/2 + o_p(h^2)$. The $o_p(h^2)$ term is eliminated by the assumption that $Lh^5 = O(1)$. □

LEMMA A.7. *With all the assumptions above, we have that*
$$\widehat{\mathcal{V}}(x_1, x_2, \Delta) = \widehat{\mathcal{V}}_0(x_1, x_2, \Delta) + O_p(L^{-1}h^{-1/2}).$$

PROOF. Notice that
$$\widehat{\mathcal{V}}(x_j, x_l, \Delta) = \widehat{\mathcal{V}}_0(x_j, x_l, \Delta)$$
$$+ \left[\sum_r \{\overline{Y}_{r\cdot j} - \Psi_r(x_j)\}c_r(x_l, \Delta)\right.$$
(16)
$$+ \sum_r \{\overline{Y}_{r\cdot l} - \Psi_r(x_l)\}c_r(x_j, \Delta)$$
$$\left. + \sum_r \{\overline{Y}_{r\cdot j} - \Psi_r(x_j)\}\{\overline{Y}_{r\cdot l} - \Psi_r(x_l)\}a_r(\Delta)\right]a^{-1}(\Delta),$$
$$c_r(x_1, \Delta) = L^{-1}\int\int \{Y(s_1, x_1) - \Psi_r(x_1)\}K_h\{\Delta - (s_1 - s_2)\}N_2(ds_1, ds_2).$$

Using the expression above, it is easy to see that $E\{c_r(x_1, \Delta)\} = 0$, and calculations as in Lemma A.3 show that
$$\text{var}\{c_r(x_1, \Delta)\} = L^{-2}\int\int\int\int K_h\{\Delta - (s_1 - s_2)\}K_h\{\Delta - (s_3 - s_4)\}$$
$$\times [\mathcal{V}\{x_1, x_1, (s_1 - s_3)\} + I(s_1 = s_3)\sigma_\varepsilon^2]$$
$$\times E\{N_2(ds_1, ds_2)N_2(ds_3, ds_4)\}$$
$$= O(\nu^2 L^{-1}h^{-1}).$$

On the other hand, $\overline{Y}_{r\cdot j} - \Psi_r(x_j) = \frac{1}{N_r}\int\{Y_r(s,x_j) - \Psi_r(s,x_j)\}N(ds)$. It is easy to see that $E\{\overline{Y}_{r\cdot j} - \Psi_r(x_j)\} = 0$, and that

$$\text{var}[N_r\{\overline{Y}_{r\cdot j} - \Psi_r(x_j)\}]$$
$$= \int\int [\mathcal{V}\{x_j, x_j, (s_1 - s_2)\} + I(s_1 = s_2)\sigma_\varepsilon^2]$$
$$\quad \times \{\nu^2 g(s_1/L)g(s_2/L)\,ds_1 ds_2 + \nu g(s_1/L)\varepsilon_{s_1}(ds_2)\,ds_1\}$$
$$= \nu^2 L^2 \int \mathcal{V}(x_j, x_j, Lu)f_1(u)\,du + \nu L \int \{\mathcal{V}(x_j, x_j, 0) + \sigma_\varepsilon^2\}g(s_1)\,ds_1$$
$$= \nu^2 L f_1(0)\int \mathcal{V}(x_j, x_j, u)\,du + \nu L\{\mathcal{V}(x_j, x_j, 0) + \sigma_\varepsilon^2\} + o(L).$$

By properties of Poisson processes, we have $N_r/(\nu L) \to 1$ a.s. Therefore, we have $\overline{Y}_{r\cdot j} - \Psi_r(x_j) = O_p(L^{-1/2})$, $c_r(x_1, \Delta) = O_p(L^{-1/2}h^{-1/2})$. By Lemma A.2, $a_r(\Delta) = O_p(1)$. Therefore, $\widehat{\mathcal{V}}(x_1, x_2, \Delta) - \widehat{\mathcal{V}}_0(x_1, x_2, \Delta) = O_p(L^{-1}h^{-1/2})$, completing the proof. $\square$

PROOF OF THEOREM 1. This is a direct result from Lemmas A.6 and A.7. $\square$

PROOF OF THEOREM 2. For a fixed $\Delta \neq 0$, when $h \leq |\Delta|$, we have $\widetilde{\mathcal{V}}(x_1, x_2, \Delta) = \{\widehat{\mathcal{V}}(x_1, x_2, \Delta) + \widehat{\mathcal{V}}(x_1, x_2, -\Delta)\}/2$. This equation is true automatically for $\Delta = 0$. Therefore, the asymptotic distribution of $\widetilde{\mathcal{V}}(\Delta)$ is the same as that of $\{\widehat{\mathcal{V}}(\Delta) + \widehat{\mathcal{V}}(-\Delta)\}/2$, for any fixed $\Delta$.

For $\Delta_1 \neq \pm\Delta_2$, by Theorem 1, $\{\widehat{\mathcal{V}}(\Delta_1), \widehat{\mathcal{V}}(-\Delta_1)\}^T$ and $\{\widehat{\mathcal{V}}(\Delta_2), \widehat{\mathcal{V}}(-\Delta_2)\}^T$ are asymptotically independent, and the joint asymptotic normality of the four vectors can be established. Therefore $\widetilde{\mathcal{V}}(\Delta_1)$ and $\widetilde{\mathcal{V}}(\Delta_2)$ are jointly asymptotic normal and asymptotically independent. It suffices to show that $\Omega(\Delta)$ is the asymptotic covariance matrix of $\widetilde{\mathcal{V}}(\Delta)$.

For $\Delta \neq 0$, apply the delta method to the joint asymptotic distribution of $\widehat{\mathcal{V}}(\Delta)$ and $\widehat{\mathcal{V}}(-\Delta)$; the following gives the asymptotic covariance between $\widetilde{\mathcal{V}}(x_1, x_2, \Delta)$ and $\widetilde{\mathcal{V}}(x_3, x_4, \Delta)$:

$$(1/4)(RLh)^{-1}\{\nu^2 f_1(0)\}^{-1}R_K$$
$$\times \{\mathcal{M}(x_1, x_2, x_3, x_4, \Delta, \Delta, 0) + \mathcal{M}(x_1, x_2, x_3, x_4, -\Delta, -\Delta, 0)$$
$$+ 2\mathcal{M}(x_1, x_2, x_3, x_4, \Delta, -\Delta, -\Delta) + 2I(x_2 = x_4)\sigma_\varepsilon^2 \mathcal{V}(x_1, x_3, 0)$$
$$+ 2I(x_1 = x_3)\sigma_\varepsilon^2 \mathcal{V}(x_2, x_4, 0) + 2I(x_1 = x_3, x_2 = x_4)\sigma_\varepsilon^4$$
$$+ 2I(x_2 = x_3)\sigma_\varepsilon^2 \mathcal{V}(x_1, x_4, 0) + 2I(x_1 = x_4)\sigma_\varepsilon^2 \mathcal{V}(x_2, x_3, 0)$$
$$+ 2I(x_1 = x_4, x_2 = x_3)\sigma_\varepsilon^4\}.$$

Note that $\mathcal{M}(x_1, x_2, x_3, x_4, -\Delta, -\Delta, 0) = \mathcal{M}(x_1, x_2, x_3, x_4, \Delta, \Delta, 0)$ by the symmetry in the definition of $\mathcal{M}(x_1, x_2, x_3, x_4, u, v, w)$. Next, for $\Delta = 0$, we have $\widetilde{\mathcal{V}}(x_1, x_2, 0) = \widehat{\mathcal{V}}(x_1, x_2, 0)$, and the asymptotic covariance between $\widetilde{\mathcal{V}}(x_1, x_2, 0)$ and $\widetilde{\mathcal{V}}(x_3, x_4, 0)$ is given in Theorem 1. The proof is completed. $\square$

PROOF OF COROLLARY 1. The result follows from Theorem 2 and the delta method. To see this, note that, with the separable structure in (3), we have $\mathcal{V}(x_1, x_2, \Delta) = G(x_1, x_2)\rho(\Delta)$ and $\mathcal{V}^{(0,0,2)}(x_1, x_2, \Delta) = G(x_1, x_2)\rho^{(2)}(\Delta)$. By the delta method, the asymptotic mean of $\widehat{\rho}(\Delta)$ is

$$\frac{\sum_{x_1 \in \mathcal{X}} \sum_{x_2 \leq x_1} \{\mathcal{V}(x_1, x_2, \Delta) + \sigma_K^2 \mathcal{V}^{(0,0,2)}(x_1, x_2, \Delta) h^2/2 + o_p(h^2)\}}{\sum_{x_1 \in \mathcal{X}} \sum_{x_2 \leq x_1} \{G(x_1, x_2) + \sigma_K^2 G(x_1, x_2) \rho^{(2)}(0) h^2/2 + o_p(h^2)\}}$$
$$= \{\rho(\Delta) + \sigma_K^2 \rho^{(2)}(\Delta) h^2/2 + o_p(h^2)\}/\{1 + \sigma_K^2 \rho^{(2)}(0) h^2/2 + o_p(h^2)\}$$
$$= \{\rho(\Delta) + \sigma_K^2 \rho^{(2)}(\Delta) h^2/2 + o_p(h^2)\} * \{1 - \sigma_K^2 \rho^{(2)}(0) h^2/2 + o_p(h^2)\}$$
$$= \rho(\Delta) + \{\rho^{(2)}(\Delta) - \rho(\Delta)\rho^{(2)}(0)\}\sigma_K^2 h^2/2 + o_p(h^2).$$

The asymptotic variance of $\widehat{\rho}(\Delta)$ also follows from the delta method. $\square$

**Acknowledgments.** The work of Raymond Carroll occurred during a visit to the Center of Excellence for Mathematics and Statistics of Complex Systems at the Australian National University, whose support is gratefully acknowledged. We also thank the three anonymous referees for providing valuable and constructive suggestions.

Y. Li
Department of Statistics
University of Georgia
Athens, Georgia 30602
USA
E-mail: yehuali@uga.edu

M. Hong
N. D. Turner
J. R. Lupton
Faculty of Nutrition
Texas A&M University
College Station, Texas 77843-2253
USA
E-mail: meeyoung_hong@yahoo.com
n-turner@tamu.edu
jlupton@tamu.edu

N. Wang
Department of Statistics
Texas A&M University
College Station, Texas 77843-3143
USA
E-mail: nwang@stat.tamu.edu

R. J. Carroll
Department of Statistics
Texas A&M University
College Station, Texas 77843-3143
USA
and
Faculty of Nutrition
Texas A&M University
College Station, Texas 77843-2253
USA
E-mail: carroll@stat.tamu.edu